\documentclass[11pt]{article}

\usepackage{epsfig,amsmath,latexsym}
\usepackage{amsfonts,amssymb}

\newtheorem{theorem}{Theorem}[section]

\newtheorem{lemma}[theorem]{Lemma}

\def\slfrac#1#2{\hbox{\kern.1em %
 \raise.5ex\hbox{\the\scriptfont0 #1}\kern-.11em %
 /\kern-.15em\lower.25ex\hbox{\the\scriptfont0 #2}}}

\newcommand{\pf}{\noindent{\bf Proof.~}}
\newcommand{\beq}{\begin{eqnarray}}
\newcommand{\eeq}{\end{eqnarray}}
\newcommand{\beql}[1]{\begin{eqnarray}\label{#1}}
\newcommand{\beqs}{\begin{eqnarray*}}
\newcommand{\eeqs}{\end{eqnarray*}}
\newcommand{\eqn}[1]{(\ref{#1})}

\newcommand{\cc}{{\mathbb C}}
\newcommand{\rr}{{\mathbb R}}
\newcommand{\zz}{{\mathbb Z}}

\newcommand{\RR}{{\mathbb R}}
\newcommand{\ZZ}{{\mathbb Z}}

\newcommand{\ft}{{\mathfrak t}}

\newcommand{\br}{{\mathbf r}}

\newcommand{\bg}{{\mathbf g}}

\newcommand{\bo}{{\mathbf 1}}

\newcommand{\bv}{{\mathbf v}}

\newcommand{\bz}{{\mathbf z}}

\newcommand{\bA}{{\mathbf A}}
\newcommand{\bB}{{\mathbf B}}
\newcommand{\bC}{{\mathbf C}}
\newcommand{\bD}{{\mathbf D}}
\newcommand{\bG}{{\mathbf G}}
\newcommand{\bI}{{\mathbf I}}
\newcommand{\bJ}{{\mathbf J}}
\newcommand{\bM}{{\mathbf M}}
\newcommand{\bP}{{\mathbf P}}
\newcommand{\bQ}{{\mathbf Q}}
\newcommand{\bS}{{\mathbf S}}
\newcommand{\bT}{{\mathbf T}}
\newcommand{\bU}{{\mathbf U}}
\newcommand{\bV}{{\mathbf V}}
\newcommand{\bW}{{\mathbf W}}
\newcommand{\bX}{{\mathbf X}}
\newcommand{\bY}{{\mathbf Y}}
\newcommand{\bZ}{{\mathbf Z}}

\newcommand{\sA}{{\mathcal A}}

\newcommand{\sD}{{\mathcal D}}
\newcommand{\sG}{{\mathcal G}}

\newcommand{\sL}{{\mathcal L}}

\newcommand{\sS}{{\mathcal S}}

\newcommand{\bsq}{\vrule height .9ex width .8ex depth -.1ex}

\makeatletter
\def\@sect#1#2#3#4#5#6[#7]#8{\ifnum #2>\c@secnumdepth
     \def\@svsec{}\else
     \refstepcounter{#1}\edef\@svsec{\csname the#1\endcsname.\hskip .75em }\fi
     \@tempskipa #5\relax
      \ifdim \@tempskipa>\z@
        \begingroup #6\relax
          \@hangfrom{\hskip #3\relax\@svsec}{\interlinepenalty \@M #8\par}%
        \endgroup
       \csname #1mark\endcsname{#7}\addcontentsline
         {toc}{#1}{\ifnum #2>\c@secnumdepth \else
                      \protect\numberline{\csname the#1\endcsname}\fi
                    #7}\else
        \def\@svsechd{#6\hskip #3\@svsec #8\csname #1mark\endcsname
                      {#7}\addcontentsline
                           {toc}{#1}{\ifnum #2>\c@secnumdepth \else
                             \protect\numberline{\csname the#1\endcsname}\fi
                       #7}}\fi
     \@xsect{#5}}
\def\@begintheorem#1#2{\it \trivlist \item[\hskip \labelsep{\bf #1\ #2.}]}

\def\plain{plain}\ifx\fmtname\plain\csname fi\endcsname
     
     \input docstrip
     \preamble

     Do not distribute the stripped version of this file.
     The checksum in the header refers to the documented version.

     \endpreamble
     \generateFile{here.sty}{t}{\from{here.doc}{}}
     \endinput
\fi
\ifcat a\noexpand @\let\next\relax\else\def\next{%
    \documentstyle[here,doc]{article}\MakePercentIgnore}\fi\next
\ifx\@Hxfloat\@Hundef\else\expandafter\endinput\fi
\let\@Hxfloat\@xfloat
\def\@xfloat#1[{\@ifnextchar{H}{\@HHfloat{#1}[}{\@Hxfloat{#1}[}}
\def\@HHfloat#1[H]{%
\expandafter\let\csname end#1\endcsname\end@Hfloat
\vskip\intextsep\vbox\bgroup\def\@captype{#1}\parindent\z@
\ignorespaces}
\def\end@Hfloat{\egroup\vskip \intextsep}

\makeatother

\catcode`\@=11

\catcode`@=12

\begin{document}


\begin{center}
{\Large {\bf Apollonian Circle Packings: Geometry and Group Theory}}\\
{\Large {\bf II. Super-Apollonian Group and Integral Packings}}\\
\vspace{1.5\baselineskip}
{\em Ronald L. Graham}\\
Department  of Computer Science and Engineering\\
University of California at San Diego, 
La Jolla, CA 92093-0114 \\
\vspace*{1.5\baselineskip}

{\em Jeffrey C. Lagarias} \\
Department of Mathematics \\
University of Michigan, 
Ann Arbor, MI 48109--1109 \\
\vspace*{1.5\baselineskip}

{\em Colin L. Mallows} \\
Avaya Labs, Basking Ridge, NJ 07920 \\
\vspace*{1.5\baselineskip}

{\em Allan R. Wilks} \\
AT\&T Labs, Florham Park, NJ 07932-0971 \\
\vspace*{1.5\baselineskip}

{\em Catherine H. Yan}
\footnote{
Partially supported by NSF grants DMS-0070574, DMS-0245526 and a Sloan
Fellowship. This author is also affiliated with Dalian University of
Technology, China.}\\
Department of Mathematics \\
Texas A\&M University,
College Station, TX 77843-3368\\
\vspace*{1.5\baselineskip}
(March 10, 2005 version) \\
\vspace*{1.5\baselineskip}
{\bf ABSTRACT}
\end{center}
Apollonian circle packings arise by repeatedly filling the interstices
between four mutually tangent circles  with further tangent circles.
Such packings can be described in terms of the Descartes configurations
they contain, where a  
Descartes configuration is a set of four mutually tangent
circles in the Riemann sphere, having disjoint interiors.
Part I showed there exists a discrete group, the Apollonian group,
acting on a parameter space of 
(ordered, oriented) Descartes configurations, 
such that the Descartes configurations in a packing formed an
orbit under the action of this group. It is observed there
exist infinitely many types of integral
Apollonian  packings in which all circles had integer curvatures,
with the integral structure being related to the integral
nature of the Apollonian group. 
Here we consider the
action of a larger discrete group, the super-Apollonian
group, also having an integral structure,
whose orbits describe the Descartes quadruples of a geometric
object we call a 
super-packing. The circles in a super-packing never cross each other
but are nested to an arbitrary depth.
Certain Apollonian packings and super-packings
are strongly integral in the sense that the curvatures of all
circles are integral and the curvature$\times$centers of all
circles are integral. We show that (up to scale)
there are exactly $8$ different
(geometric) strongly integral super-packings, and that each
contains a copy of every integral Apollonian circle packing
(also up to scale). We show  that the super-Apollonian group
has finite volume in the group of all automorphisms of the
parameter  space of 
Descartes configurations, which is isomorphic
to the Lorentz group $O(3, 1)$.

\vspace*{1.5\baselineskip}

\setcounter{page}{1}

%
%
%
\setlength{\baselineskip}{1.0\baselineskip}
\section{Introduction}
\setcounter{equation}{0}

\begin{figure}[htbp]
\centerline{\epsfxsize=5.0in \epsfbox{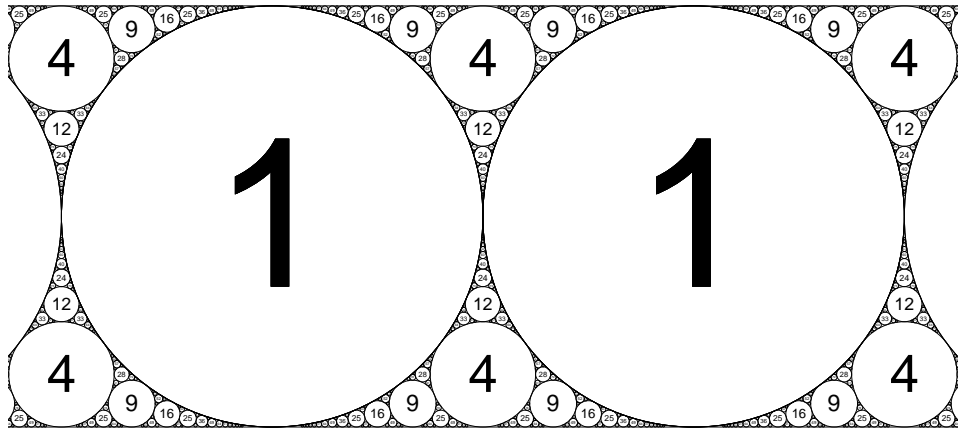}}~\label{fig1}
\caption{The integer  Apollonian packing $(0,0, 1,1)$}
\end{figure}

Apollonian circle packings are arrangements of
tangent circles that arise by repeatedly filling the interstices
between four mutually tangent circles  with further tangent circles.
A set of four mutually tangent circles is called a
Descartes configuration. 
Part I studied Apollonian circle packings in terms of
the set of Descartes configurations that they contain.
It is observed that there exist  Apollonian circle packings 
that have a very strong integral structure, with all
circles in the packing having integer curvatures,
and rational centers, such that curvature$\times$center
is an integer vector. We termed these {\em strongly
integral} Apollonian circle packings. An example is 
the $(0,0,1,1)$ packing pictured in Figure~\ref{fig1},
with the two circles of radius $1$ touching at the origin, and with
two straight lines parallel to the $x$-axis.

%
%

Part I gave an explanation for the existence of such
integral structures. This uses  a coordinate
system for describing all (ordered, oriented)
Descartes configurations $\sD$ in terms of their 
curvatures and centers, which forms a $4 \times 3$
{\em curvature-center coordinate}
matrix $\bM_{\sD}$, and a more detailed coordinate system,
{\em augmented curvature-center coordinates},
using $4\times 4$ matrices $\bW_{\sD}$. 
The strongly integral property of a single Descartes configuration
is encoded in the integrality of the matrix $\bM_{\sD}$.
The set of all (geometric) Descartes configurations in an 
Apollonian packing can be described as a single  orbit of
a certain discrete group $\sA$ of $4 \times 4$ 
integer matrices; algebraically there are $48$  orbits
of ordered, oriented Descartes configurations giving rise to 
the same geometric packing, which  correspond to the 
48 possible ways of ordering the circles and totally orienting the
configuration.
The integrality of the members
of $\sA$ is the source of the strong integrality
of some Apollonian circle packings. 
As a consequence of this group action, if  a single
Descartes configuration in the packing is strongly
integral, then they all are; hence  
every individual circle in the
packing is strongly integral. 

There are infinitely many distinct integral Apollonian
circle packings. Two more  of them are pictured in
Figures 2 and 3,  the $(-1, 2,2,3)$-packing
and the $(-6, 11, 14, 15)$-packing, respectively.
 Any integral packing
can be moved by a Euclidean motion so as to be
strongly integral, as will follow from results in
this paper.

%
%

\begin{figure}[htbp]
\centerline{\epsfxsize=3.0in \epsfbox{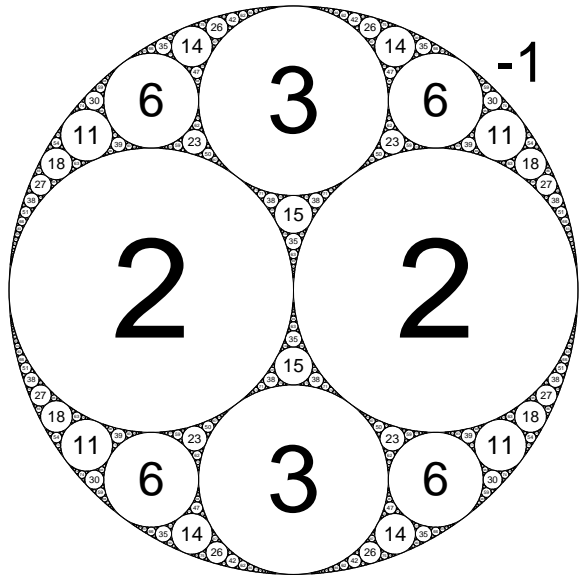}}
\caption{The integer  Apollonian packing $(-1,2,2,3)$}~\label{fig2}
\end{figure}
%
%
%
\begin{figure}[htbp]
\centerline{\epsfxsize=3.0in \epsfbox{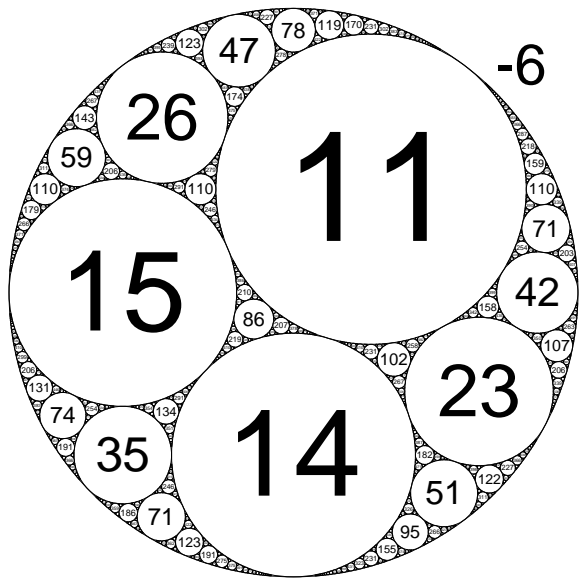}}
\caption{The integer  Apollonian packing $(-6,10,11,14)$}~\label{fig3}
\end{figure}

Part I introduced some further integral actions on
Descartes configurations, involving a ``duality''
operator $\bD$ leading to a dual Apollonian group $\sA^{\perp}$.
Combining this group with the Apollonian group
leads to a large group of integer $4 \times 4$ matrices,
the {\em super-Apollonian group} $\sA^{S}$, which 
acts on the set of all (ordered, oriented) Descartes
configurations. Part I showed that the super-Apollonian
group is a hyperbolic Coxeter group. 
It defined an {\em Apollonian super-packing}
to be an orbit  $\sA^{S}[\sD]$ of a single
Descartes configuration $\sD$ under this group. 
Such a super-packing is called  {\em integral}
if the initial Descartes configuration has integer
curvatures, and is  
{\em strongly integral} if the curvature-center coordinates
 $\bM_{\sD}$ of the initial Descartes configuration
are integral. These properties hold for all Descartes
configurations in the packing if they hold for one.

\begin{figure}[htbp]
\centerline{\epsfxsize=6.0in \epsfbox{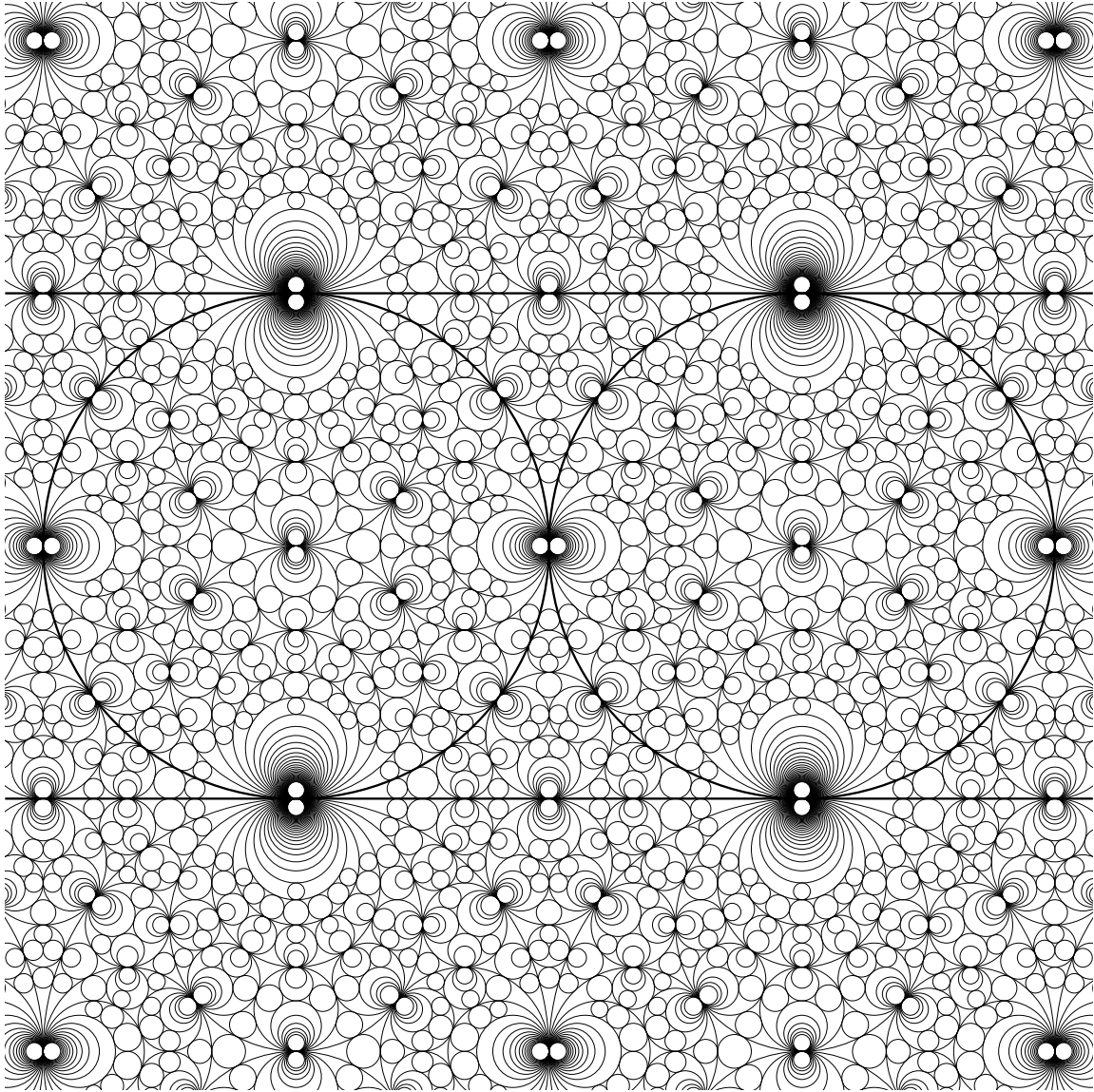}}
\caption{An initial part of the  $(0,0,1,1)$ super-packing
(square of sidelength $4.4$)}~\label{fig4}
\end{figure}

 In this paper we study the geometric 
structure  of super-packings,
and their integrality properties.
A  striking geometric fact is
 that the circles in a super-packing do not
overlap, as shown in \S3. Figure~\ref{fig4}  shows
circles of curvature at most $65$ in the 
super-packing generated from the $(0,0,1,1)$ 
configuration in Figure~\ref{fig1} above.
(The generating Descartes configuration is indicated
with slightly darker lines.)
Here the horizontal and vertical scales of the
figure are roughly from $-2.2 \le x \le 2.2$
and $-2.2 \le y \le 2.2$. This picture is 
representative of  the
whole super-packing; we show in \S6 that this
super-packing is periodic under the 
two-dimensional lattice generated
by the shifts
$x \to x+2$ and $y \to y+2$.  

%
%
%

We cannot picture the whole super-packing because
the circles in it are dense in the plane; they nest inside
each other to an arbitrary depth. One can show that
the set of points lying in the interior of an infinite
nested sequence of such circles has full Lebesgue
measure. However there is interesting structure
visible if we turn up the scale of the circles.
In \S6 in Figure~\ref{fig8} we show all circles of curvature at
most 200 inside the super-packing above, inside
the unit square $0 \le x \le 1,~ 0 \le y \le 1$.

This paper is concerned with integrality 
properties of some super-packings, and we study successively
stronger integrality properties. In \S4 we require the
curvatures to be integral, in \S5 and \S6 both
the curvatures and curvature$\times$center must be 
integral, and in \S8 the augmented curvature-center
coordinates of all Descartes configurations in the 
super-packing must be integral. At this last level we
show that there are exactly $14$ such super-packings, viewed
as rigid geometric objects.
In \S2 we give a more detailed description of our results.

The general framework of this paper 
was developed by  the second author (JCL), who also did much of
the writing. This paper is an extensive
revision  of a preprint written in 2000, and adds new
results in \S3 and  \S6.

\paragraph {Acknowledgments.} The authors are grateful for
helpful comments from Andrew Odlyzko, Eric Rains,
Jim Reeds and  Neil Sloane during this work. They
also thank the  reviewer for incisive comments. \\

%
%
%
\section{Summary of Main Results}
\setcounter{equation}{0}

We consider Apollonian super-packings.
Analogously to part I,
an   Apollonian  super-packing may be considered as either
an geometric object or an algebraic object, as follows. \\

(i) [Geometric] A {\em geometric Apollonian super-packing} 
is a point set on the Riemann sphere
$\hat{\cc} = \RR^2 \cup \{ \infty\}$, consisting of 
all the circles in four orbits
of a certain group $G_{\sA^{S}}(\sD)$ of
M\"{o}bius transformations
inside  the conformal group M\"{o}b(2) acting
on the four circles $\{C_1, C_2, C_3, C_4\}$ in
a given Descartes configuration $\sD$.  The 
group $G_{\sA^{S}}(\sD)$ depends on $\sD$. \\

(ii) [Algebraic] An {\em (algebraic) Apollonian super-packing} is a
 set of ordered, oriented Descartes configurations,
given by $48$ orbits of the
super-Apollonian group $\sA^{S}[\sD]$. The augmented 
curvature-center coordinates of
its elements are $\sA^{S}\bW_{\sD}:=\{ \bU W_{\sD} : \bU \in \sA^{S}\}$. \\

A geometric super-packing  can be 
described  in terms of its unordered, unoriented
Descartes configurations. From this viewpoint,
 each geometric Apollonian
super-packing corresponds to  $48$ different algebraic 
super-packings; there are 
$24$ choices of ordering the four circles 
and $2$ choices of total orientation of the configuration. 
We can  consider geometric  super-packings as unions
of a countable  number of Apollonian packings. 
This leads to interesting questions concerning the way these
Apollonian packings are embedded inside the geometric super-packing.
We note that as  a point set, the geometric super-packing is
invariant under the group action $G_{\sA^{S}}(\sD)$. 
However it is not a closed set,   and its closure is
the entire Riemann sphere $\RR^2 \cup \{ \infty\}$.

A large part of the paper considers 
integrality properties of curvatures and centers
of some super-packings. These questions are mainly
studied using algebraic super-packings, although
we also consider questions concerning 
the associated geometric super-packing, such as
its group of symmetries under Euclidean motions.

In \S3  Theorem~\ref{th31} shows that each  geometric 
super-packing is 
a packing in the sense that the circles in it do not cross
each other transversally, as mentioned above.
Theorem~\ref{th32n} specifies certain subcollections of 
geometric super-packings which are genuine packings
in the sense that the interiors of the circles do not
overlap.

In \S4 we study integer super-packings,
in which all circles have integer curvatures.
Integer super-packings are classified up to
Euclidean motions by a 
single invariant, their {\em divisor} $g$,
which is the greatest common divisor of the
curvatures in any Descartes configuration in
the packing. Theorem~\ref{th41n} shows  that for
each $g \ge 1$ there is a unique such
integral super-packing, up to a Euclidean
motion. We also show in Theorem~\ref{th43n} 
that for each geometric Apollonian circle 
packing that is integral, there exists a
Euclidean motion taking it to one that is
strongly integral.

In \S5 we study strongly integral super-packings,
which are those whose curvatures are integral and whose
curvature$\times$center is also integral.
Strongly-integral super-packings are geometrically rigid:
Theorem~\ref{orbits} shows that for each integer $g \ge 1$ there are
exactly $8$ strongly integral geometric super-packings
which have divisor $g$. Here we do not allow the packings
to be moved by Euclidean motions.

\begin{figure}[htbp]
\centerline{\epsfxsize=3.0in \epsfbox{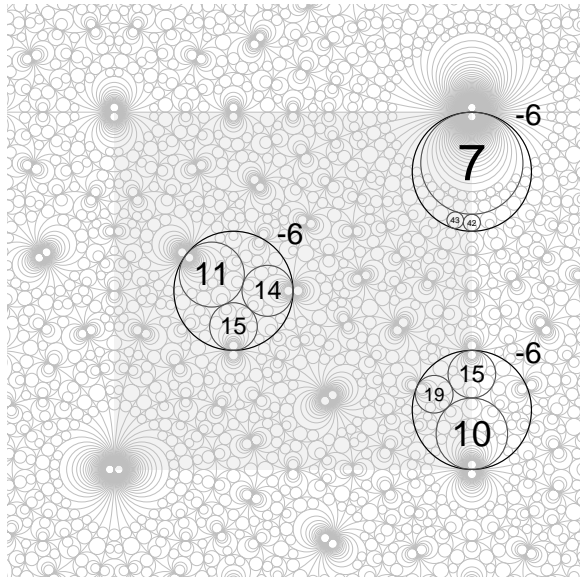}}
\caption{Integer 
Apollonian packings with bounding circle of curvature $6$.}~\label{fig5}
\end{figure}

\begin{figure}[htbp]
\centerline{\epsfxsize=3.0in \epsfbox{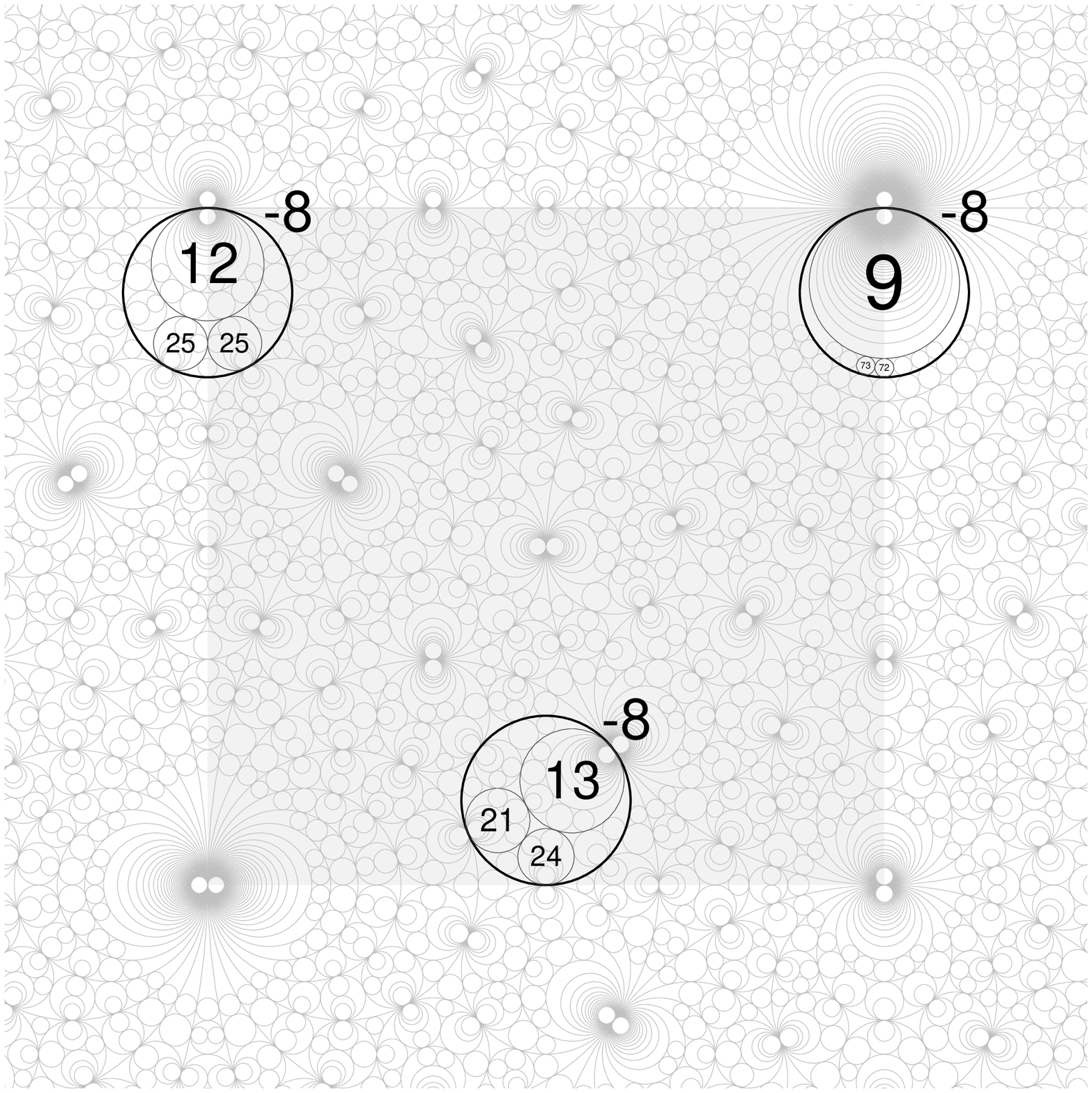}}
\caption{Integer Apollonian packings with bounding circle of curvature $8$.
}~\label{fig6}
\end{figure}

\begin{figure}[htbp]
\centerline{\epsfxsize=3.0in \epsfbox{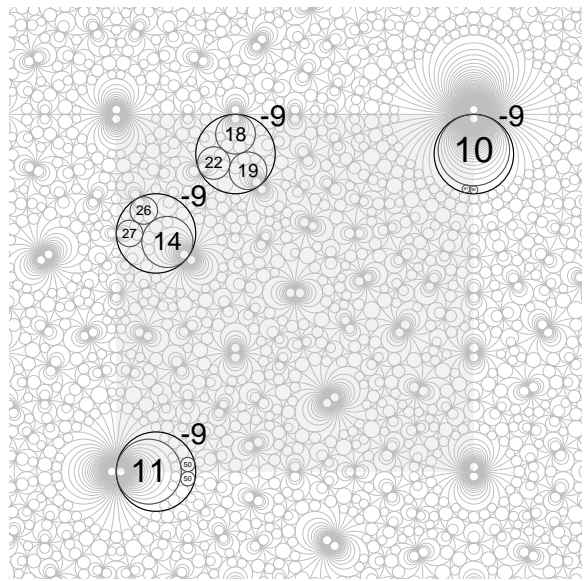}}
\caption{Integer Apollonian packings with bounding circle of curvature $9$.
}~\label{fig7}
\end{figure}

In \S6 we study the relations between integer Apollonian packings
and strongly integral super-packings. 
Without loss of generality we restrict to
primitive integer super-packings, those with divisor $1$.
Theorem~\ref{th61} shows
 that each of the $8$ kinds of these has a large group of
internal symmetries, forming a crystallographic group.
For convenience we fix one of them and call it the
{\em standard strongly integral super-packing}; results proved for
it have analogues for the other seven. Theorem~\ref{th62}
shows that each primitive integral Descartes configuration
(except for the $(0,0,1,1)$ 
configuration) occurs in this packing with the center of its largest
circle being contained in the closed unit square
$0 \le x \le 1, ~ 0 \le y \le 1$, 
and the location of this circle center is unique.  
Theorem~\ref{th63} deduces  that the geometric 
standard strongly integral  super-packing contains a
unique copy of each primitive integral Apollonian packing, 
except for  the
$(0,0,1,1)$ packing, 
having the property that 
the center of its bounding circle lies inside the
closed unit square.
Figures~\ref{fig5}, \ref{fig6} and \ref{fig7} picture
the locations of all primitive integer Apollonian
packing with bounding circle of curvatures 6, 8 and
9, respectively. The unit square is indicated by 
slightly darker shading in the figures. 
Note that in Figure~\ref{fig5} the three Apollonian packings
are generated by Descartes configurations with
curvature vectors $(-6, 7, 42, 43),~ (-6, 10, 15,19)$
and $(-6,11, 14, 15)$; these are root quadruples in 
the sense of \cite[Sec. 4]{GLMWY2}.

%
%
%
 
%
%
%
 
%
%
%

In \S7 Theorem~\ref{th51} shows  that the super-Apollonian
group is a finite index normal subgroup of the 
group $Aut(Q_{D}, \ZZ)$ of integral automorphs of 
the Descartes quadratic form. The latter group can
be identified with an index $2$ subgroup of the
integer
Lorentz group $O(3, 1, \zz)$, and this identification
allows us to identify 
the super-Apollonian group with  a 
particular normal subgroup $\tilde{\sA}^{S}$
of index $96$ in  $O(3, 1, \zz)$, defined after
 Theorem~\ref{th51}.

In  \S8 we study  super-packings all of whose 
Descartes configurations are {\em super-integral} in the
sense that their augmented curvature-center coordinate
matrices $\bW_{\sD}$ are integer matrices. Theorem~\ref{th81n} 
shows  there
are exactly $14$ geometric super-packings of this kind.

The last section \S9 makes a few concluding remarks.

The Appendix considers minimal conditions to guarantee
that a Descartes configuration
is strongly integral.  Theorem~\ref{th101} shows
that a configuration is strongly integral if (and only
if) three of its four circles are strongly integral.

%
%
%
\section{Geometric Apollonian Super-Packings}
\setcounter{equation}{0}

In this section we consider properties of 
geometric Apollonian circle packings. We view such
a super-packing as a point set on the Riemann sphere
$\hat{\cc} = \cc \cup \{\infty\}.$  We first
note  that it is not a closed set.
It  is not hard to show that its  closure is
the whole Riemann sphere, 
Each geometric Apollonian packing  has  a group invariance
property under a certain group of M\"{o}bius transformations
which depends on the super-packing, which is 
the group generated by the countable set of  
groups of M\"{o}bius transformations that leave 
invariant some Apollonian packing contained in
the super-packing. 

Our object is  prove the following
``packing'' property of geometric super-packings.
%
%
%
%

\begin{theorem}~\label{th31} 
A  geometric Apollonian super-packing is a circle
packing in the weak sense that no two circles belonging
to it cross each other transversally. Circles in the geometric 
super-packing may be nested, or tangent to each other.
\end{theorem}

Before giving the proof, we describe the 
nature of the geometric packing in terms of 
nesting of circles. We view the packing 
$\sA^{S}[\sD_0]$ as generated from an initial
(positively oriented) Descartes configuration $\sD_0$,
by multiplication by a finite set of generators of
the super-Apollonian group.
Each  circle in
the super-packing has a  well-defined 
{\em nesting depth} $d$ 
(relative to the generating configuration $\sD_0$)
which counts the number of circles in the packing which
include $C$ in their interior. Here the notion of
``interior'' is defined with respect to  the initial Descartes
configuration $\sD_0$.
The Apollonian group
generators move ``horizontally'', leaving constant
the nesting depth of any circles they produce. 
The dual Apollonian group generators move ``vertically'',
by reflecting three of the circles in a configuration
into the interior of the fourth circle, they increase
the nesting depth by one. We show there is a unique
``normal form'' word of minimal length in the generators
that produces a Descartes
configuration $\sD$ containing $C$.
The nesting depth of $C$ 
is exactly equal to the number $d$ of generators
of $\sA^{\perp}$ that appear in this normal form word.
The circles at nesting depth $0$
are those circles in the Apollonian packing generated by
$\sD_0$. Each circle $C$ at nesting depth $k \ge 1$ contains
a unique Apollonian packing, consisting  of it
plus all circles at depth $k+1$ nested inside it.

\paragraph{Proof of Theorem~\ref{th31}.} 
We view the geometric super-packing on the Riemann sphere
$\hat{\cc}= \cc \cup \{ \infty \},$ so the initial 
Descartes configuration $\sD_0$ 
consists of four circles on the sphere. In this case each
circle defines a spherical cap. We choose an ordering
and orientation of $\sD_0$ (this does not affect the geometry),
requiring that  $\sD_0$ have  positive (total)  orientation.

Let the super-packing be $\sA^{S}[\sD_0] = \sA^{S} \bW_{\sD_0}.$
We  consider the effect of the super-Apollonian group generators
acting on the left on
the matrix $\bW_{\sD_0}$.
The Descartes
configurations in the super-packing are given by 
$\sD = \bU_m \bU_{m-1} \cdots \bU_1[\sD_0]$
with
each $\bU_k \in \{ \bS_1, \bS_2, \bS_3, \bS_4,
\bS_1^{\perp}, \bS_2^{\perp}, \bS_3^{\perp}, \bS_4^{\perp}\}$.
We consider words in the group measured by their {\em length} $m$.
The stage $m$ circles will consist of all new circles added
using products of $m$ generators. 
We may without loss of generality restrict to 
{\em normal form} words, which are those that satisfy
the two conditions:

(i) If $\bU_k = \bS_i$, then $\bU_{k+1} \ne \bS_i$
and $\bU_{k+1} \ne \bS_j^{\perp}$ with $j \ne i$.

(ii) If $\bU_k = \bS_i^{\perp}$, then  $\bU_{k+1} \ne \bS_i^{\perp}$. 

Equivalently, looking backwards, if $\bU_{k+1} = \bS_i^{\perp}$
then either $\bU_k= \bS_i$ or else $\bU_{k}=\bS_j^{\perp}$ 
for some $j \ne i$. A word may be put in
normal form by canceling adjacent equal generators, since
all $\bS_i^2 = (\bS_i^{\perp})^2= \bI$, and by moving towards
the right\footnote{That is, moving it towards
the beginning of the word.}
in the word as far as possible any generator $\bS_i^{\perp}$,
using the property that it commutes with all $\bS_j$ with  $j \ne i$. 
These operations eventually put a word in normal form,
without increasing its length. The operations do not change
the Descartes configuration $\sD$ it represents.

We prove the theorem by induction on the
number of symbols $m$ in a normal form word,
which we call the stage of the induction.
The induction hypotheses at stage $m$ are  as follows.
Here we  let $C_i$ refer to
the circle at row $i$ of the associated
(ordered, oriented) Descartes configuration.

(1) Each normal form word of length $m$ produces
either one or three new circles, according as
$\bU_{m}=\bS_i$, where it is the circle
$C_i$ or $\bU_{m}= \bS_j^{\perp}$, where it is
the three circles $C_i$ with $i \ne j$. 

(2) The nesting depth of any new circle produced
at this stage is equal to the number of occurrences of a letter
in $\sA^{\perp}$ in its generating normal form
word.

(3) Each new circle produced has empty interior,
when it first appears. 

In particular hypothesis (3) implies that all circles
produced at stage $m$ have disjoint interiors, and 
that each such circle  
contains no circles from  earlier stages in 
its interior. 
If the induction is proved, then 
hypothesis  (3) guarantees  that the nesting depth of a
circle is well-defined when it is first produced,
because no new circle will ever include it in its
interior. Hypothesis  (3) also  guarantees  that all
circles produced by the end of stage $m+1$ do
not cross. As a consequence  no two circles in
the packing cross, using  property (3) 
applied at that level $m$ which is 
the greater of the levels of the two circles.
Thus the theorem will follow.

The base case $m=0$ of the
induction is immediate, consisting of the
initial Descartes configuration $\sD_0$.
We now show the induction step for $m+1$, given $m$. 
We are given  a Descartes configuration
$\sD' = \bU_{m+1} [\sD] = \bU_{m+1} \bU_m \cdots \bU_1[\sD_0]$
with a normal form word. 

To establish hypothesis (1) for $m+1$, suppose first
that  $\bU_{m+1} = \bS_i^{\perp}.$
We assert that  the $i$-th circle $C_i$ of $\sD$ was a new
circle produced at stage $m$. For either $\bU_m = \bS_i$,
in which case it was the unique new circle in $\sD$
by induction hypothesis (1) at stage $m$, or $\bU_m = \bS_j^{\perp}$
with $j \ne i$, in which case it was one of three new circles
produced at stage $m$. By induction hypothesis (2) the
circle $C_i$ has empty interior at stage $m$. The
three circles $\{C_j^{'}: j \ne i\}$ in 
the new configuration $\sD'$  are contained
in the interior of  $C_i$ are therefore new circles.
They do not cross, being part of a Descartes configuration.
Thus hypothesis (1) holds for $m+1$ in this case. 

Suppose next that $\bU_{m+1} = \bS_i$,
so the possible new circle is $C_i^{'}$.
If 
$\sD_k:=\bU_k \bU_{k-1} \cdots \bU_1[\sD_0]$ 
is  the maximal 
length subword
such that $\bU_k= \bS_j^{\perp}$ for some $j$, 
with $1 \le k \le m$ then 
$\sD'$ belongs to the Apollonian packing generated by
$\sD_k$, since all subsequent generators belong to the
Apollonian group. If no such $k$ exists, then
$\sD'$ is in the Apollonian packing generated
by the original Descartes configuration $\sD_0$.
For $k \ge 1$, this  Apollonian packing is
entirely contained  in
the interior of a 
bounding circle $C= C_j$ first produced
at stage $k-1$. At that time the interior of $C$ was empty,
by induction hypotheses (2) and (3) at stage $k-1$.
The only Descartes configurations that ever can enter
the interior of the circle $C$, must do so by
a reflection in $C$ and these   are exactly
those normal form words
starting with initial segment $\bU_k \bU_{k-1} \cdots \bU_1$.
(This follows from uniqueness of a circle when it is
created, hypothesis (3) applied at stage $k$.)
The words contain this initial segment at the same depth,
with all subsequent letters in $\sA$, fill out an Apollonian
packing at this depth.
 In particular each such normal form word produces one
new circle in this Apollonian packing.
Recall that  all the circles
in the Apollonian packing inside the bounding circle $C$
have disjoint interiors (Theorem 4.1 of part I). 
These circles all have the same depth, and (1) holds in
this case. Normal form words that have another subsequent
generator in $\sA^{\perp}$ confine the resulting 
Descartes configuration to the
inside  of a  single circle 
in the Apollonian packing $\sA[\sD_k]$ already produced,
and all longer words with this prefix are entirely
contained inside this circle. In particular,
they do not coincide with  $C_i^{'}$.
the new circle produced by the normal form word corresponding
to $\sD'$.
It follows that the circle $C_i^{'}$ is new,
so hypothesis (1) holds for $m+1$ in this case.

There remains the case where $\bU_{m+1} = \bS_i$
and all previous $\bU_k \in \sA$. Then 
$\bU_{m+1} \cdots \bU_1[\sD_0]$ belongs
to the Apollonian packing generated by
$\sD_0$, and is at depth $0$. Here $C_i^{'}$
is new since the generation of an 
Apollonian packing creates one new circle
at each step, and any word that contains
an element of $\sA^{\perp}$ moves the
Descartes configuration inside an older
circle in this packing, from which it cannot
escape. Thus $C_i^{'}$ is a new circle in
this case, and hypothesis (1) holds.

Hypothesis (2) holds for $m+1$ in the case $\bU_m =\bS_i^{\perp}$
because the three new circles produced have nesting
depth one greater than the nesting depth of $C_i$ at
the previous level, to which induction hypothesis (2) applied.
It also holds in the remaining case $\bU_m=\bS_i$ because the 
argument above showed that the nesting 
depth did not increase in this case.

Hypothesis (3) holds if $\bU_{m+1} =\bS_i^{\perp}$. Now 
circle $C_i$ was first created only at stage $m$,
and the only other possible sequences leading to
a Descartes configuration including 
this circle 
must start from $\sD_m$ and use a generator
$\bU_{m+1}= \bS_j$
with  $j \ne i$.  In all other cases the resulting
Descartes configuration includes no circle inside
$C_i$, so the interiors of the three new circles 
produced are empty at the end of stage $m+1$.

In the remaining case  $\bU_{m+1} =\bS_i$, we have
already observed that the new circle $C_i^{'}$ produced,
is disjoint from all other circles in the Apollonian
packing  $\sA[\sD_k]$ created by level $m$, 
which are necessarily contained  in the
bounding circle $C$ of $\sD_k$. 
As mentioned above, all Descartes configurations
containing a circle inside $C$ must have an
initial segment of their generating word, giving
the unique normal form word that first generates $C$,
at stage $k \le m$. Other depth $m+1$ words with 
this initial segment, and with all subsequent $\bU_j$ drawn 
from the Apollonian group generators 
produce new circles in the
Apollonian packing $\sA[\sD_k]$, disjoint from $C_i^{'}$.
Any  normal  form word
at depth $m+1$ with this initial segment, which 
contain  some generator $\bS_j^{\perp}$ after-wards,
produce a Descartes configuration contained inside a circle
of the Apollonian packing $\sA[\sD_k]$
different from $C$. 
It follows that the
interior of $C_i^{'}$ is empty at the end of stage
$m+1$, as required.

This completes the induction step and the proof.
$~~~\bsq$ \\

Super-packings have  some properties that
are genuine packing properties. The proof of Theorem~\ref{th31}
established in hypothesis  (3) shows that the finite set of circles at
stage $m$ of the construction, starting from a generating
Descartes configuration $\sD$,  all had disjoint interiors.
It also gives  the following stronger result. 

%
%
%
%

\begin{theorem}~\label{th32n} 
For a geometric Apollonian super-packing given with
a generating Descartes configuration $\sD$, 
and each $k \ge 1$,  
the set of all circles having  nesting depth exactly $k$ 
with respect to $\sD$ have pairwise  disjoint interiors.
These circles can be viewed as  
forming an infinite collection
of Apollonian packings, each missing one circle;
the missing circle is a bounding circle at depth $k-1$.
\end{theorem}

\paragraph{Proof.} 
The proof of Theorem~\ref{th31}
shows that the nesting depth of a
circle is well-defined.  Circles at nesting depth $k$ have
disjoint interiors, since no two of them cross, and the
only way for two of them to have an interior point in
common is for one to be nested inside the other, which 
would violate the nesting ordering. 

Given a normal form word $\bU := \bU_m \bU_{m-1} \cdots \bU_1$
with $\bU_m = \bS_i^{\perp}$, and containing exactly
$k$ elements of 
$\sA^{\perp} = 
\langle \bS_1^{\perp},\bS_2^{\perp}, \bS_3^{\perp}, \bS_4^{\perp} \rangle, $
the set of all normal form words having $\bV$ as prefix and all
other letters $\bU_j \in \sA$ for $j > m$ produce all the
circles in an Apollonian packing $\sA[\sD_m]$, 
with $\sD_m = \bU[\sD_0]$.
All these circles are  at depth $k$ 
except for  the outer circle of $\sD_m$, which 
is its $i$-th circle.
Enumerating all possible such $\bU$ as prefixes  represents
the set of nesting depth $k$ circles as a collection of
Apollonian packings, each excluding one circle when $k \ge 1$.
$~~~\bsq$ \\

\paragraph{Remarks.}
(1) The proof of Theorem~\ref{th31} is a geometric analogue of the
presentation for the super-Apollonian group
proved in part I, \cite[Theorem 6.1]{GLMWY11}.

(2) The analogous result to Theorem~\ref{th31}
fails 
to hold in all dimensions $n \ge 4$, as explained
in part III. The ``nesting'' property
of the dual Apollonian configurations resulting
from ``vertical'' moves still exists
and works in all dimensions.
However the ``horizontal'' motions  moving spheres  around
in Apollonian packings produces  spheres that cross
in dimensions $n \ge 4$, see \cite[Lemma 4.1]{GLMWY13}.

(3) We cannot easily visualize a geometric super-packing as
a completed object because the circles in it are dense
in the plane. We can however picture a partial version
of it that pictures all circles of size above a given
threshold, in some finite region of the packing.
The integral super-packings we are most interested
in have a periodic lattice of symmetries (see Theorem~\ref{th61}),
so it suffices to examine a finite region of the packing.
Figure~\ref{fig4} in \S1 and Figure~\ref{fig8} in \S6
exhibit part of a super-packing. 

(4) Every circle $C$ in a geometric super-packing $\sA^{S}[\sD]$ 
has associated
to it a unique Apollonian packing of which it is the bounding
circle. If  it is a depth $k$ circle (relative  to the
starting configuration $\sD$), then this Apollonian packing 
consists of all depth $k+1$ circles contained in the interior of
$C$.

%
%

\section{Integral Super-Packings}
\setcounter{equation}{0} 

An Apollonian super-packing is {\em integral} if it
contains one (and hence all) Descartes configuration
whose circles have integer curvatures.

An invariant of an integral super-packing is its
{\em divisor} $g$, which is the greatest common divisor
of the curvatures of the circles in any Descartes
configuration in the super-packing.
The quantity  $g$ is well-defined   independent of
the Descartes configuration chosen, using the
relation $\bM_{\sD'} = \bU \bM_{\sD}$ between two
such configurations, where $\bU \in \sA^{S}$. Since $\bU$
is an integer matrix with determinant $\pm 1$, it preserves
the greatest common divisor of each column of the integer
matrix $\bM_{\sD}$. Here  the first column
encodes the curvatures of the circle. 

%
%
%
%

\begin{theorem}~\label{th41n} 
For each integer $g \ge 1$ there exists an
integral Apollonian super-packing with 
divisor $g$. The associated geometric super-packing
is unique, up to a Euclidean motion.
\end{theorem}

As an immediate corollary of this result, the
integral super-packing with divisor $g$ contains
at least one  copy of every integral Descartes configuration
having divisor $g$. For each such a Descartes configuration
generates an integral super-packing with divisor $g$,
so the corollary follows by the uniqueness assertion.

We defer the proof of 
Theorem \ref{th41n} to the end of the section.
It is based on a reduction theory which finds inside
any such integral super-packing a Descartes configuration
having particularly simple curvatures. 
%
%
%
%

\begin{theorem}~\label{reduction}
Let $\sD$ be an integral Descartes configuration, with
divisor $g := \gcd(b_1, b_2, b_3, b_4)$.
Then the Apollonian super-packing $\sA^{S} [\sD]$
generated by $\sD$ contains a Descartes configuration having
curvature vector a permutation of either $(0,0, g, g)$
or $(0,0, -g,-g)$, with the former case occurring if
$b_1+b_2+b_3+b_4 > 0$ and the latter case if $b_1+b_2+b_3+b_4 < 0$. 
\end{theorem}

\noindent {\bf Proof.} 
Since the super-Apollonian group $\sA^S$ preserves the (total) orientation 
of Descartes configurations, it is sufficient to show that for a positively 
oriented integral Descartes configuration $\sD$ with curvatures
$(b_1, b_2, b_3, b_4)$, $b_1+b_2+b_3+b_4 > 0$, there exists $\bU \in \sA^S$ 
and a permutation matrix $\bP_\sigma$ such that
$$\bP_\sigma \bU(b_1, b_2, b_3, b_4)^T=(0,0,g,g)^T.$$

We measure
the {\em size} of the curvature vector $\bv = (b_1, b_2, b_3, b_4)^T$
of a  Descartes configuration  by
$$ size(\bv ):= \bo^T \bv = b_1 + b_2 + b_3 + b_4.$$
We claim that for positively oriented integral Descartes configurations
with greatest common divisor $g$, we have
$$size(\bv)\geq 2g,$$
and equality holds if and only if $\bv$ is a permutation of $(0,0,g,g)$.
If all curvatures are nonnegative this is clear, since at most two
can be zero, and the other two are positive integers.
Now in any Descartes configuration  at most one circle can have
negative curvature, call it $b_1=-a$  ($a \in \ZZ_+$), in which case  
it encloses the other three.
Each of these three enclosed circles has a larger curvature in
absolute value than
the bounding circle, so  $b_i \geq a + 1$ for $i=2,3,4$. Thus
$size(\bv) \geq -a+ 3(a + 1) \geq 2a+3 > 2g,$
which proves the claim.

We give a reduction procedure which
chooses  matrices in $ \sA^{S}$
to reduce the size and show that the procedure halts only at a
vector of form $(0, 0, g,g)$, up to a permutation.
To specify it, we observe that for the curvature vector 
$\bv=(b_1, b_2, b_3, b_4)^T$ of 
any integral
Descartes configuration  with
$b_1 \leq b_2\leq b_3\leq b_4$, we have
\beql{605aa}
size(\bS_4 \bv) = \bo^T \bS_4\bv \leq \bo^T \bv = size(\bv).
\eeq
and equality holds
if and only if $b_1b_2+b_2b_3+b_3b_1=0$. To see this,
we have $\bS_4\bv=(b_1, b_2, b_3, b_4')^T$ where
$$b_4'=2(b_1+b_2+b_3)-b_4 =b_1+b_2+b_3-2\sqrt{b_1b_2+b_2b_3+b_3b_1}.$$
Thus
$$\bo^T \bS_4 \bv-\bo^T \bv=b_4'-b_4=-4\sqrt{b_1b_2+b_2b_3+b_3b_1} \leq 0.$$
Equality can hold if and only if $b_1b_2+b_2b_3+b_3b_1=0$, which proves the
observation. Note that $b_1 \leq 0$ when the equality in \eqn{605aa} holds.
Also note that $g=\gcd(b_1, b_2, b_3, b_4)$ is an invariant under the 
action of $\sA^S$. 

Starting with any  positively oriented integral Descartes configuration with
curvature vector $(b_1, b_2, b_3, b_4)^T$, 
where $b_i$ is the largest number, we apply
$\bS_i \in \sA$.
By \eqn{605aa},  the $size(\bv)$ decreases but cannot be negative,
so after a finite
series of $\bS_i$ we arrive at a positively oriented 
integral Descartes configuration with curvatures 
$\bv'= (b_1', b_2', b_3', b_4')^T$,
where $\gcd(b_1', b_2', b_3', b_4')=g$ and
the smallest curvature, say $b_1'$, satisfies $b_1' \leq 0$, and 
the size of $\bv'$ can not be reduced by the action of the Apollonian group
$\sA$.
Call this the basic reduction step. Note that the basic reduction step
involves only matrices in the Apollonian group and therefore
moves around inside a single Apollonian packing.

If $b_1' =0$ then necessarily $b_2'=0$,
whence the curvature vector is $(0,0,b_3',b_3')^T$,
and by $g=\gcd(b_1',b_2', b_3', b_4')$, we have 
$b_3' = g$ and the reduction halts.
If $b_1'<0$, applying $\bS_1^T$, we get a new Descartes configuration with
$\bv''= (-b_1', b_2'+2b_1', b_3'+2b_1', b_4'+2b_1')^T$,
which is positively oriented and lies in a new Apollonian packing
and has $$size(\bv'') = size(\bv) + 4 b_1' < size(\bv).$$
Thus the size strictly
decreases and is non-negative. Now we may re-apply the basic reduction step.
Continuing in this way we get strict decrease of $size(\bv)$ at each
step, with the only possible halting step being the smallest curvature 
equals $0$. Since the size of the curvature vector 
is bounded below and decreases by at least one at each
step, the procedure terminates at $(0,0,g,g)$, up to a permutation.~~~$\bsq$

\paragraph{Proof of Theorem~\ref{th41n}.}
For existence, the  super-packing generated by
a Descartes configuration with curvature vector $(0,0, g,g)$,
which is a homothetically scaled version of the configuration
$(0,0,1,1)$ pictured in Figure 1, is necessarily integral 
with divisor $g$.

For uniqueness,
Theorem~\ref{reduction} shows that any two geometric integral Apollonian
super-packings with divisor $g$ each contain a Descartes configuration
whose curvatures are $(0,0,g,g)$ up to permutation and orientation. 
Now  it is true for any two such Descartes configurations with identical
curvature vectors are congruent,
i.e. one is obtainable from the other by a Euclidean motion.
This is obvious by inspection for the $(0,0,g,g)$ Descartes
configuration, which necessarily
consists  of two touching circles of radius $\frac{1}{g}$
and two parallel lines.

Now the Euclidean motion that takes one
Descartes configuration to the other, also takes the super-packing
generated by the first configuration  to
the one generated by the  other, because the super-packing is defined by the
action of the super-Apollonian group on the left on the
Descartes configuration $\bW_{\sD}$, and this commutes with
the Euclidean motion acting as a M\"{o}bius transformation
on the right. This establishes uniqueness.
~~~$\bsq$

We can use the freedom of
a Euclidean motion allowed in Theorem~\ref{th41n}
to make an internal Apollonian super-packing strongly integral.

\begin{theorem}~\label{th43n}
For each integral geometric Apollonian super-packing 
there is a Euclidean
motion that takes it to a strongly integral geometric
Apollonian super-packing.
\end{theorem}

\paragraph{Proof.}
Using Theorem~\ref{reduction}  
each integral geometric Apollonian packing contains 
a Descartes configuration $\sD$ with curvatures $(0,0, g,g)$; 
note that for a  geometric packing the order and orientation of the
Descartes configuration do not matter.  The curvature vector
determines the Descartes configuration up to 
congruence. We can now find a strongly integral
Descartes configuration $\sD'$ with this curvature vector.
For  $g=1$ such a configuration is given explicitly by
\eqn{N601} below, and for larger $g$ we obtain a strongly
integral configuration from it using the homothety 
$(x, y) \mapsto \frac{1}{g}(x, y)$. There exists a Euclidean
motion that maps $\sD$ to $\sD'$, since they are congruent
configurations. This motion maps the super-packing $\sA^{S}[\sD]$
to the super-packing $\sA^{S}[\sD']$,  which is strongly integral.
~~~$\bsq$

%
%
\section{Strongly Integral Super-packings} 
\setcounter{equation}{0}

A Descartes configuration $\sD$ is {\em strongly integral}
if its associated $4 \times 3$
curvature center-coordinate matrix 
$\bM_{\sD}$ is an integer matrix; this property is
independent of ordering or orientation of the
Descartes configuration.
A super-packing is called {\em strongly integral} if
it contains one (and hence all) Descartes configurations 
having 
this property. Since a strongly integral super-packing
is integral, it has a divisor $g$ as an invariant.

Our main object in this section is to classify strongly
integral super-packings, as follows.

%
%
%
%
\begin{theorem}~\label{orbits}
(1) For each $g \ge 1$ there  are exactly $8$ 
different geometric Apollonian
super-packings that are strongly integral
and have divisor $g$.

(2) The set of all ordered, oriented Descartes configurations
that are strongly integral and have a given divisor $g$
fill exactly 384 orbits of the super-Apollonian group.
\end{theorem}

This theorem classifies  these super-packings as rigid objects, not
allowed to be moved by Euclidean motions. 
To prove this result we derive a normal form for
a ``super-root quadruple'' in a super-packing of the kind
above, as follows. 

%
%
%
\begin{theorem}\label{normal_form} 
Given a strongly integral Apollonian super-packing $\sA^s[\sD_0]$ with the 
 divisor  $g \geq 1$, 
there exists a unique ``reduced'' Descartes configuration 
$\sD \in \sA^S[\sD_0]$ 
whose curvature-center coordinate matrix $\bM=\bM_\sD$ is of the form
$\bA_{m,n}[g]$ or $\bB_{m,n}[g]$ for  $m, n \in \{0, 1\}$,  
up to a permutation of rows,  where
\beq \label{ABmn} 
\bA_{m,n}[g]=\pm \left[ \begin{array}{ccc} 
     0 & 0 & 1 \\
     0 & 0 & -1 \\
     g & m & n \\
     g & m-2 & n \end{array}  \right]
\qquad  
\bB_{m,n}[g]=\pm \left[ \begin{array}{ccc} 
     0 & 1 & 0 \\
     0 & -1 & 0 \\
     g & m & n \\
     g & m & n-2 \end{array}  \right],
\eeq
and the sign is determined by the orientation of $\sD_0$. 
\end{theorem}
\pf For a strongly integral Apollonian super-packing $\sA^S[\sD_0]$ with the 
divisor  $g$, 
by Theorem \ref{reduction} there exists a
strongly integral Descartes configuration $\sD \in \sA_S[\sD_0]$ with
curvatures $\pm(0,0,g,g)$. The two 
straight lines in $\sD$ must be parallel to either $x$-axis or $y$-axis.
It follows that the $4\times 3$ curvature-center coordinate matrix $\bM_\sD$
is of the form $\bA_{m,n}[g]$, or $\bB_{m,n}[g]$, for some $m, n \in \ZZ$, 
up to a permutation of rows. 

We now reduce $m, n$ to take the values $0, 1$ using the
following identities, which are easy to check.
\begin{eqnarray*}
\bP_{(34)}\bS_3 \bA_{m,n}[g]=\bA_{m-2, n}[g],  & \qquad&
\bP_{(34)}\bS_3 \bB_{m,n}[g]=\bB_{m, n-2}[g], \\
\bP_{(34)}\bS_4 \bA_{m,n}[g]=\bA_{m+2, n}[g], & \qquad &
\bP_{(34)}\bS_4 \bB_{m,n}[g]=\bB_{m, n+2}[g], \\
\bP_{(12)}\bS_1^T \bA_{m,n}[g]=\bA_{m,n+2}[g], & \qquad&
\bP_{(12)}\bS_1^T \bB_{m,n}[g]=\bB_{m+2,n}[g], \\
\bP_{(12)}\bS_2^T \bA_{m, n}[g]=\bA_{m, n-2}[g],  & \qquad &
\bP_{(12)}\bS_2^T \bB_{m, n}[g]=\bB_{m-2, n}[g].
\end{eqnarray*}
where $\bP_{(ij)}$ is the permutation matrix that exchanges $i$ and $j$. 
Also note that  $\bP_\sigma \bS_i=\bS_{\sigma(i)} \bP_\sigma$,
$\bP_\sigma \bS_i^T =\bS_{\sigma(i)}^T \bP_\sigma$.  Hence there is
a series of group operations in $\sA^S$  which takes $\bM_\sD$ to 
a permutation of $\bA_{m,n}[g]$ or $\bB_{m,n}[g]$ with $m, n \in \{0, 1\}$. 
This proves the existence of the ``reduced'' Descartes configuration 
in the Apollonian super-packing $\sA^S[\sD_0]$. 

To prove the uniqueness, it suffices to show that the 
$24 \times 8 \times 2=384$ 
Descartes configurations whose curvature-center coordinate matrices
are 
$$
\{ \bP \bA_{m,n}[g], \bP \bB_{m,n}[g] \ | \ \bP \in \text{Perm}_4, m, n \in \{0, 1\}
\}
$$ 
are in different Apollonian super-packings. 
(There are two signs for each of $\bA_{m,n}[g]$ and $\bB_{m,n}[g]$. )
In what follows we let $\tilde{\bA}_{m,n}[g]$, and $\tilde{\bB}_{m,n}[g]$
denote the unique $4 \times 4$ augmented curvature-center coordinate
matrices extending $\bA_{m,n}[g]$ and $\bB_{m,n}[g]$, respectively;
uniqueness holds  by Theorem 3.1 of part I.

Note that each $\bS_i$ and $\bS^T_i$ preserves the (total) orientation 
of the Descartes configuration, as well as the parity of every 
element of $\bM_\sD$. 
First we show that the matrices 
$\bA_{m,n}[g], \bB_{m,n}[g], (m,n \in \{0,1\}$ 
are in distinct orbits of $\sA^S \times \text{Perm}_4$.
To see this, for each integral vector $\bv \in \zz^4$, 
let $\kappa(\bv)$ be the
number of even terms in $\bv$, and
for any strongly integral Descartes configuration $\sD$ 
let 
$$
\kappa(\bM_\sD)=(\kappa(\bv_1), \kappa(\bv_2), \kappa(\bv_3)),
$$
where 
$\bv_1, \bv_2, \bv_3$ are the column vectors of $\bM_\sD$.
Then $\kappa(\bM_\sD)$ is invariant under the action of 
$\sA^S \times \text{Perm}_4$.  For $m, n \in \{0, 1 \}$,
$\kappa(\bA_{m,n}[g]), \kappa(\bB_{m,n}[g])$ are all distinct except
$\kappa(\bA_{1, 0}[g])=\kappa(\bB_{0, 1}[g])=(*, 2, 2)$, where
$*$ is $2$ if $g$ is odd, and $4$ if $g$ is even. However, 
$\bA_{1,0}[g]$ and $\bB_{0,1}[g]$ can not be equivalent under the action of 
$\sA^S \times \text{Perm}_4$. Arguing by contradiction, 
assume that there exists a
matrix $\bU \in \sA^S \times \text{Perm}_4$ such that
$\bU \bA_{1, 0}[g]=\bB_{0, 1}[g]$.  This relation lifts
to augmented curvature-center coordinates:  
$\bU \tilde{\bA}_{1,0}[g]=\tilde{\bB}_{0,1}[g]$. 
It follows that 
$\bU=(\tilde \bA_{1,0}[g])^{-1} \tilde \bB_{0,1}[g]$ is unique.
We can directly verify that for $m,n \in \ZZ$ 

\beq \label{aug-ABmn} 
\tilde \bA_{m,n}[g]=\left[\begin{array}{cccc}
2(n+1)/g &  0 & 0 &  1  \\
2(1-n)/g &  0 & 0 & -1 \\
(m^2+n^2-1)/g &   g & m & n \\
((m-2)^2+n^2-1)/g &  g & m-2 & n \end{array}\right],
\nonumber \\
\tilde \bB_{m,n}[g]=\left[\begin{array}{cccc}
2(m+1)/g &   0 & 1 & 0 \\
2(1-m)/g &   0 & -1 & 0 \\
(m^2+n^2-1)/g &   g & m & n \\
(m^2+(n-2)^2-1)/g  &   g & m & n-2 \end{array}\right], 
\eeq
by checking that these satisfy the identity of Theorem 3.2 of
part I necessary and sufficient 
to be augmented curvature-center coordinates.
Now it is easy to verify that 
 $\bP_{(14)}\bP_{(23)}\bD \tilde \bA_{1, 0}[g]=\tilde \bB_{0, 1}[g]$, 
where $\bD=-\bQ_D$ is defined as in \S3 of Part I \cite{GLMWY11}, 
\beq \label{matrixD}
\bD=\frac{1}{2} \left[ \begin{array}{rrrr} 
 -1 & 1 & 1 & 1  \\
 1 & -1 & 1 & 1  \\
 1 & 1 & -1 & 1 \\
 1 & 1 & 1 & -1 \end{array} \right].
\eeq 
By the uniqueness 
$\bP_{(14)}\bP_{(23)} \bD =\bU \in \sA^S \times \text{Perm}_4$, 
which is  impossible since $\sA^S \times \text{Perm}_4$ consists of 
integral matrices only, while $\bD$ has half integers.
In conclusion, $\bA_{m, n}[g], \bB_{m, n}[g]$,
($m, n \in  \{0, 1\}$) are in distinct orbits of
$\sA^S \times \text{Perm}_4$. 

The final step
is to show that for any permutation $\bP \neq \bI$, $\bP \bA_{m,n}[g]$
(resp. $\bP \bB_{m,n}[g]$) 
can not be obtained from $\bA_{m,n}[g]$ (resp. $\bB_{m,n}[g]$) 
by an action of $\sA^S$. That is, we claim: 
 if for  a permutation matrix $\bP \in \text{Perm}_4$,
there exists a matrix $\bU \in \sA^S$ such that
\[
\bU \bA_{m,n}[g] =\bP \bA_{m,n}[g], \qquad \text{or } \bU \bB_{m,n}[g] =\bP \bB_{m,n}[g],
\]
then $\bP=\bI$.

To establish the claim, 
consider  again the $4 \times 4$ ACC-coordinate matrices $\tilde \bA_{m,n}[g]$ 
and $\tilde \bB_{m,n}[g]$. From \S3.1 of Part I \cite{GLMWY11},  
for any Descartes configuration $\sD$, 
the curvature-center coordinate matrix $\bM_\sD$  
can be uniquely extended to a $4 \times 4$ ACC-coordinate matrix
$\bW_\sD$. It follows that  
if $\bU \bA_{m,n}[g] =\bP \bA_{m,n}[g]$, then
the equality holds for their $4\times 4$ ACC-coordinate matrices, i.e., 
$\bU \tilde \bA_{m,n}[g] =\bP \tilde \bA_{m,n}[g]$.  It implies  
 $\bU=\bP \in \sA^S \cap \text{Perm}_4$. 
However, comparing the size of $\bU$ and $\bP$, where the size of 
a matrix $\bU$ is defined as $f(\bU):= \bo^T \bU \bo$, we have
$f(\bP):={\bf 1}^T \bP {\bf 1}=4$ for $\bP \in \text{Perm}_4$,
and $f(\bU):={\bf 1}^T \bU {\bf 1} \geq 8$ for any $\bU \in \sA^S, 
\bU\neq \bI$,
(c.f.\S5 of Part I \cite{GLMWY11}).  Therefore 
the only possibility is $\bU=\bP=\bI$. The same argument applies to 
$\bB_{m,n}[g]$, and the claim follows. 

We conclude that a reduced Descartes
configuration of the
form \eqn{aug-ABmn} in any 
strongly integral Apollonian super-packing 
exists and is unique. ~~~\bsq

\paragraph{Proof of Theorem~\ref{orbits}.}
(1) Since there are 48 orbits of ordered, oriented
Descartes configurations corresponding to each
geometric super-packing, to show there are exactly
$8$ geometric super-packings it suffices to  show 
that the strongly integral Descartes configurations
form  $384$ orbits of the Apollonian group,
which we do below.

(2) We enumerate  the complete set of 
ordered, oriented Descartes configurations that
are strongly integral, with greatest common divisor $g$,
as follows. 
By Theorem \ref{normal_form},
any such Descartes configuration  
is equivalent under the action of $\sA^S$ to a permutation of a Descartes 
configuration whose $4\times 3$ curvature-center coordinate matrix is of 
the form $A_{m,n}$ or $B_{m,n}$, with $m,n \in \{0, 1\}$. 
The uniqueness of Theorem \ref{normal_form} asserts that the $24$ permutations
of $A_{m,n}[g]$ ($B_{m,n}[g]$) are all in distinct orbits of the 
super-Apollonian
group. Considering the two choices of orientations, we get 
$24 \times 8 \times 2=384$ orbits. $~~~\bsq$ \\

%
%

\section{Primitive Strongly Integral Super-packings}
\setcounter{equation}{0}

The strongly integral
superpackings having  a given curvature g.c.d. $g$ are each obtainable
from one  with $g=1$ by homothety.
A homothety $(x,y) \mapsto  r(x,y)$ changes all curvatures
by $\frac{1}{r}$ while leaving (curvature)$\times$(center)
unchanged. Applying the homothety with $r = \frac{1}{g}$ 
takes a strongly integral super-packing with g.c.d. $g$ to
one with g.c.d. equal to $1$.

We now  consider strongly integral super-packings
having $g=1$, which we call {\em primitive} super-packings.
Results for them carry over easily to those with divisor $g> 1$ by 
applying a homothety. Theorem~\ref{orbits} showed there are
exactly $8$ such packings. 

For convenience we will single
out a particular one of them and term it the {\em standard strongly
integral super-packing}.
We choose this to be the super-packing generated
by the ordered, oriented Descartes configuration
having 
\begin{equation}~\label{N601}
\bM_{\sD_1}= \left[
\begin{array}{rrr}
 0 & 0 & 1 \\
 0 & 0 & -1 \\
 1 & 1 & 0 \\
 1 & -1 & 0
\end{array}
\right]~~~~~~\mbox{so~that}~~~~~~~~
\bW_{\sD_1}= \left[
\begin{array}{rrrr}
2 &  0 & 0 & 1 \\
2 &  0 & 0 & -1 \\
0 &  1 & 1 & 0 \\
0 &  1 & -1 & 0
\end{array}
\right].
\end{equation}
This corresponds to a $(0,0,1,1)$ Descartes
quadruple, with the centers of the two circles
lying along the $x$-axis and the circles touching
at the origin $(0,0)$. The associated geometric
integral super-packing is the one pictured in
\S3. Results proved below for
the standard  super-packing apply generally to all eight 
primitive integral super-packings, using the Euclidean
motions mapping between them described after the
proof of Theorem~\ref{th61} below. 

We first show  that the geometric standard 
strongly-integral super-packing
has a large group of symmetries, which form  a 
crystallographic group of the plane.
%
%
%
%

\begin{theorem}~\label{th61}
The geometric standard strongly integral super-packing is
invariant under the following Euclidean motions:

(1) The lattice of translations $(x,y) \mapsto (x+2, y)$, 
and $(x,y)  \mapsto (x, y+2)$,

(2) The reflections $(x,y) \mapsto (-x, y)$, and 
$(x,y)  \mapsto (x, -y)$.

\noindent The crystallographic group generated by
these motions is  the complete set of
Euclidean motions leaving the geometric standard
strongly integral super-packing invariant.
\end{theorem}

\paragraph{Proof.}
The key fact used is that the action of the super-Apollonian
group on Descartes configurations commutes with the
action of Euclidean motions acting on Descartes
configurations as M\"{o}bius transformations. 
This was shown in part I \cite[Theorem 3.3(4)]{GLMWY11}.

(1) There is a 
Descartes configuration 
corresponding to $\sD_0$ shifted by $2$ in the 
$x$-direction and the $y$-direction; call them $\sD_0^{x}$
and $\sD_0^{y}$, respectively. These are given by the
actions of $\bS_{4}$ and  $\bS_{1}^{\perp}$, respectively.
Treating the $x$-shift first, we then have 
$$
\sA^{S}[\sD_0] = \sA^{S}[\sD_0^{x}]= \sA^{S}[ \ft_{x}(\sD_0^{'})],
$$
in which  $\ft_{x}: \bz \mapsto \bz +2$ is 
the Euclidean motion translation by $\bv=(2,0)$
as in  Appendix A of part I, and the
ordered, oriented Descartes configuration $\sD_0^{'}$ is
a permutation of $\sD_0$. 
Then the  geometric super-packing
associated to $\sA^{S}[\sD_0^{'}]$
is therefore identical with that of $\sA^{S}[\sD_0]$, and that of 
$\sA^{S}[ \ft_{x}(\sD_0^{'})]$ translates it by $(2,0)$.
Thus the geometric packing is invariant under this translation.
The argument for translation by $(0,2)$ is similar. 

(2) This geometric Descartes  configuration $\sD_0$ 
is invariant under the reflections $(x,y) \mapsto (-x, y)$, and 
$(x,y)  \mapsto (x, -y)$ viewed as M\"{o}bius transformations.
The effect of these transformations on the ordered,
oriented Descartes configuration is to permute its rows.
It follows as in (1) that the associated geometric 
Apollonian super-packings  are identical.

To see that these motions generate the full group of Euclidean
motions leaving the super-packing invariant, we observe
that the full group acts discontinuously on the plane. 
This is because the image of the $(0,0,1,1)$ configuration
is either left fixed, or else it moves a distance of at least two
in some direction, so that its circles do not overlap. 
Thus it must be contained in a crystallographic group whose
translation subgroup is given by (1) above. Now the only
possibilities are to extend the group by a subgroup of
the finite point group of motions leaving $(0,0)$ fixed (of order $8$) 
leaving
the lattice $\ZZ[(0,2), (2,0)]$ of translations invariant. 
Here (2) gives an extension of order $4$. No larger
extension occurs by observing that otherwise the image
of the $(0,0, 1,1)$ and $(-1, 2, 2, 3)$ configurations at
the origin would cross themselves.
$~~~\bsq$ \\

One can now check that 
the  $8$ primitive geometric strongly integral super-packings
given by Theorem~\ref{orbits} are obtained from
the standard strongly integral
super-packing by $8$ cosets of the Euclidean motions 
$(x, y) \mapsto (x+1, y)$,
$(x,y) \mapsto (x, y+1)$ and $(x,y) \mapsto (y, x)$
with respect to the symmetries in Theorem~\ref{th61}.
They are specified by the location 
and orientation of the $(0,0,1,1)$
configuration.

Our next result  shows that every primitive integral
Descartes configuration with no curvature zero
occurs inside the geometric standard strongly integral super-packing
in a specified location.

%
%
%
 
\begin{theorem}~\label{th62}
In the geometric standard strongly integral super-packing,
for each (unordered) primitive integral Descartes quadruple
$(a, b, c, d)$ except for $(0,0,1,1)$, 
there exists a Descartes configuration
having these curvatures, such that the
center of its largest circle lies in
the closed unit square
$\{(x,y): 0 \le x \le 1,~ 0 \le y \le 1\}$.
The location of the center of this largest circle is unique.
\end{theorem}

\paragraph{Proof.}
To establish  existence,
we first  show that a Descartes configuration
of the curvatures occurs somewhere inside 
the standard integral super-packing. This holds
because the super-packing generated by such
a configuration is an integral super-packing
with divisor $1$, which by Theorem~\ref{th41n} 
is unique up to a Euclidean motion. Thus the
standard strongly integral super-packing must
contain an isometric copy of it. Now that
we have such a configuration inside the packing,
we can use  
the translation symmetries in Theorem~\ref{th61}
to move it so that its largest circle
has center 
inside the half-open square $\{ (x, y) : -1 \le x < 1, -1 \le y < 1\}$.
If we have $-1\le x < 0$ then we apply the symmetry
$(x,y) \mapsto (-x,y)$, while if $-1\le y < 0$ we
apply $(x,y) \mapsto (x,-y)$, as necessary.

To establish uniqueness, we argue by contradiction.
If uniqueness failed,  there would exist a Euclidean motion taking
one of these Descartes configurations to the other. 
Since a single Descartes configuration generates the
entire super-packing, we conclude that the super-packing
is left invariant under this extra Euclidean motion. 
Theorem~\ref{th61} described all such automorphisms
and all of them, except the identity,  
map every point in the interior of the unit square 
square strictly outside the square. This 
contradicts the assumption that the center of
the largest circle of the first configuration
is mapped to that of the second, when
at least one of these points is strictly  inside the
unit square. In the remaining cases where both 
centers lie on the boundary, one shows 
they are must lie  the same
boundary edge and that the automorphism leaves this
edge fixed, so they are identical. 
Note that this  argument shows in passing
that the largest circle is unique, once $(0,0,1,1)$
is excluded. 
$~~~\bsq$ \\

We come now to the main result of this section,
which asserts that the 
geometric standard super-packing
contains a copy of every integral Apollonian circle
packing in a canonical way. The circles in the geometric
standard super-packing can be foliated into 
a union of geometric Apollonian packings.  

%
%
%

\begin{theorem}~\label{th63}
(1) Each circle in the standard super-packing
with center inside the 
half-open unit square
$\{(x,y): 0 \le x < 1,~ 0 \le y < 1\}$
is the exterior boundary circle of a unique 
primitive integral Apollonian circle packing contained
in the geometric standard integral super-packing.

(2) Every primitive integral Apollonian circle packing,
except for the packing  $(0,0,1,1)$,
occurs exactly once in this list. 
\end{theorem}

\paragraph{Proof.}
(1) Let the circle $C$ be given.
Recall from the proof of Theorem~\ref{th31}
that  there is a unique minimal length admissible sequence
$\bU_m \bU_{m-1} \cdots \bU_1 [\sD_0]$  
of generators of the super-Apollonian
group that yield a Descartes configuration $\sD$ containing
the given circle $C$, say in its $j$-th position.
Admissibility
requires that  $\bU_m = \bS_j$ or $\bS_i^{\perp}$ for some $i \ne j$.
Then multiplying by $\bS_j^{\perp}$ also gives an 
admissible sequence, and the Descartes configuration
$$
\sD':=\bS_j^{\perp}[\sD] =  \bS_j^{\perp} \bU_m \bU_{m-1} \cdots \bU_1 [\sD_0]
$$
consists of the circle
$C$  plus three new circles nested inside the interior of $C$.
This Descartes configuration $\sD'$ generates 
an Apollonian packing having the circle $C$ as outer
boundary, contained in the standard strongly integral super-packing. It
is unique, because if there were a second Apollonian packing
inside the bounding circle it would contain circles crossing
those in the first packing, contradicting Theorem~\ref{th31}.

(2) Recall from \cite[Sect. 3 and 4]{GLMWY2}
that inside  each integer Apollonian packing is a 
positively oriented Descartes configuration
whose absolute values of curvatures $(a, b, c, d)$ are minimal, which is
called a {\em root quadruple}. 
Theorem 4.1 of \cite{GLMWY2} showed that aside   from the root quadruple
$(0,0,1,1)$ every root quadruple is of the form 
$a< 0< b \le c \le d$. Root quadruples are characterized
by satisfying the extra condition
\beql{603root}
a+ b+c  \ge d.
\eeq
All the circles in
the resulting Apollonian packing are contained inside
a {\em bounding circle} of curvature $N=|a|$, i.e., radius $\frac{1}{N}$.
Theorem~\ref{th62} shows that for each root quadruple,
with all curvatures nonzero 
there is a matching Descartes configuration
whose largest circle has center inside the unit square
and this largest circle is unique. The Apollonian
packing contained inside this circle by  is the integral packing with
the given root quadruple, and it is unique
by the result of (1).  (In some cases, like $(-1, 2,2, 3)$
the root configuration is not unique, but the root quadruple
and the packing itself
are  always unique.) Thus every primitive integer Apollonian
circle packing, except $(0,0,1,1)$, occurs exactly once bounding
circles having center in  the closed unit square.
$~~~\bsq$ \\

The initial part of the standard  super-packing to depth 200 
inside the unit square is pictured in Figure~\ref{fig8} below.

%
%
%
\begin{figure}[htbp]
\centerline{\epsfxsize=6.0in \epsfbox{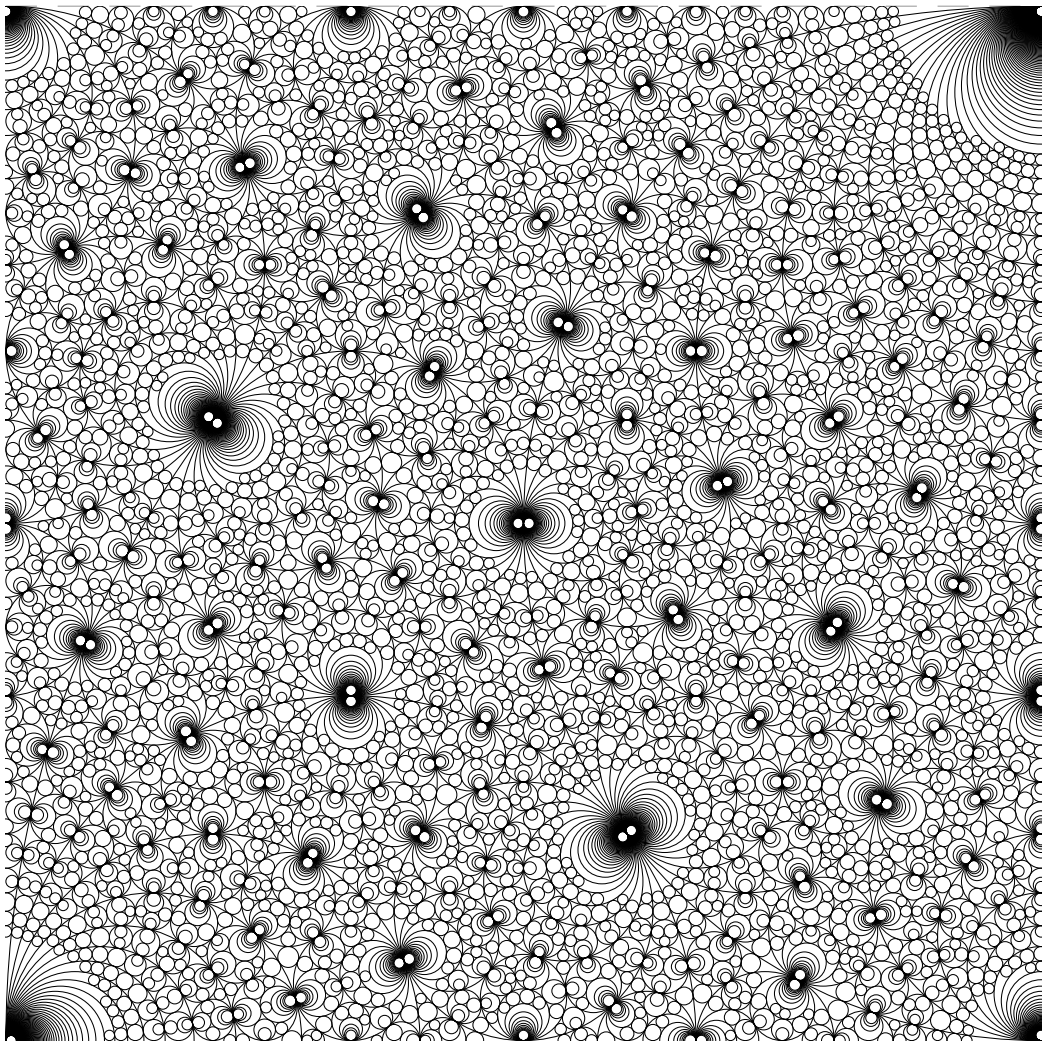}}
\caption{A ``deeper'' initial  part of a super-packing
(square of sidelength  $1$)}~\label{fig8}
\end{figure}

One can make  further 
computer experiments plotting the circles having various
curvatures restricted $(\bmod~4)$ inside the unit
square.  The results for circles
having curvatures $1~(\bmod~2)$, 
$2~(\bmod~4)$  and $0~(\bmod~4)$ and size at most 
$200$ are pictured in the following three figures.

%
%
%

\begin{figure}[htbp]
\centerline{\epsfxsize=6.0in \epsfbox{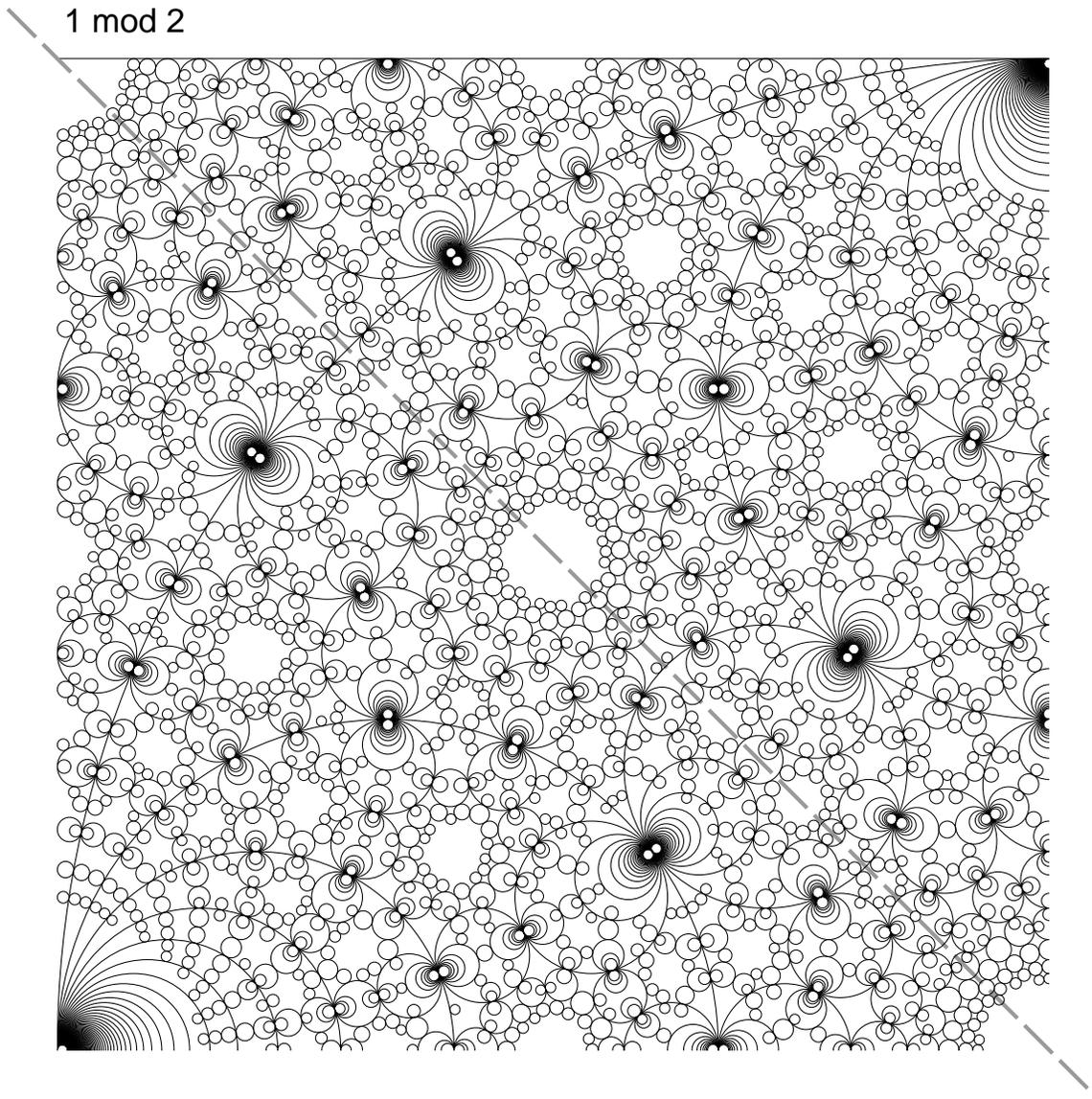}}
\caption{Circles of curvature $1~(\bmod~2)$}~\label{fig9}
\end{figure}

%
%
%

\begin{figure}[htbp]
\centerline{\epsfxsize=6.0in \epsfbox{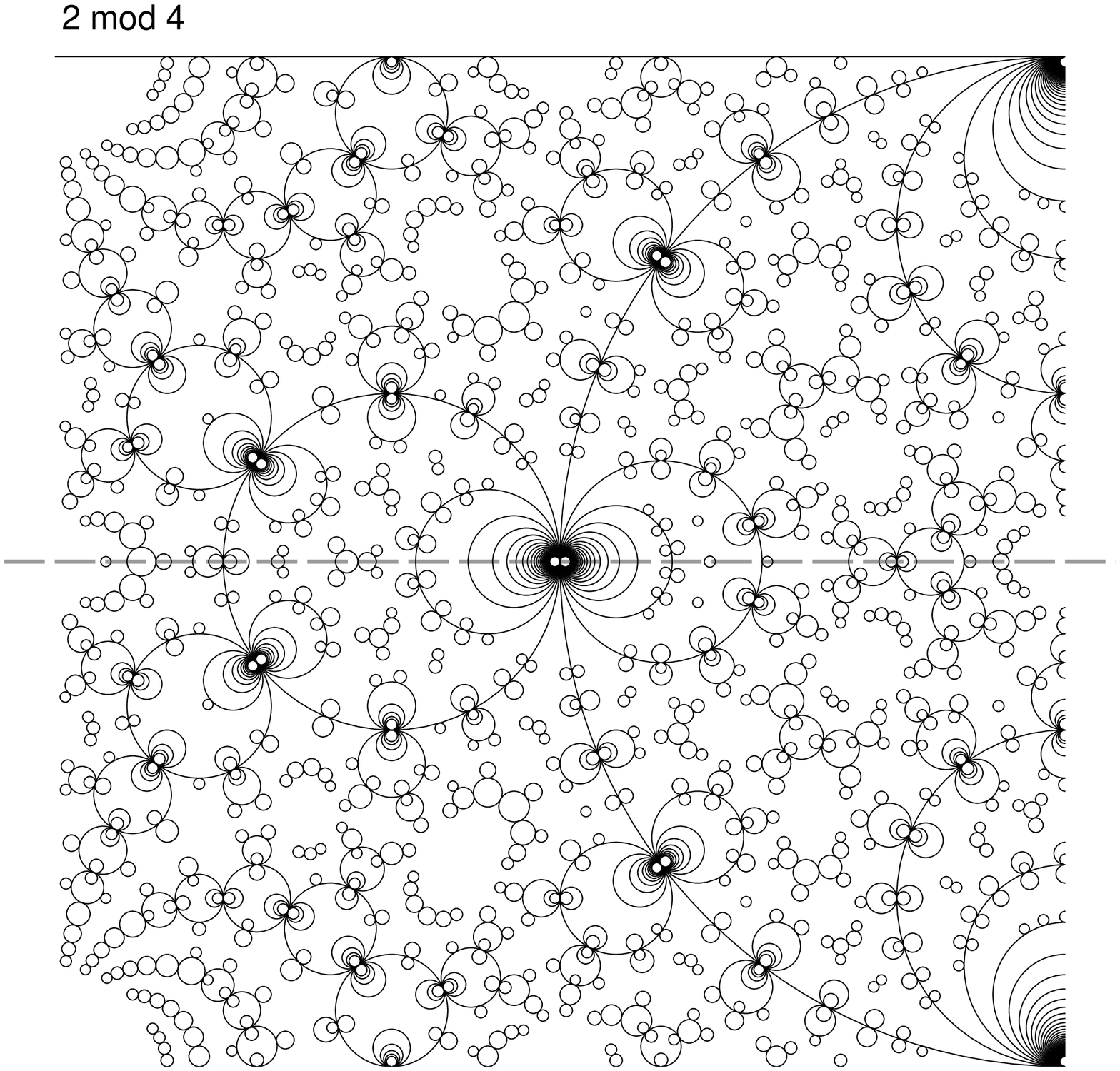}}
\caption{Circles of curvature $2~(\bmod~4)$}~\label{fig10}
\end{figure}

%
%
%

\begin{figure}[htbp]
\centerline{\epsfxsize=6.0in \epsfbox{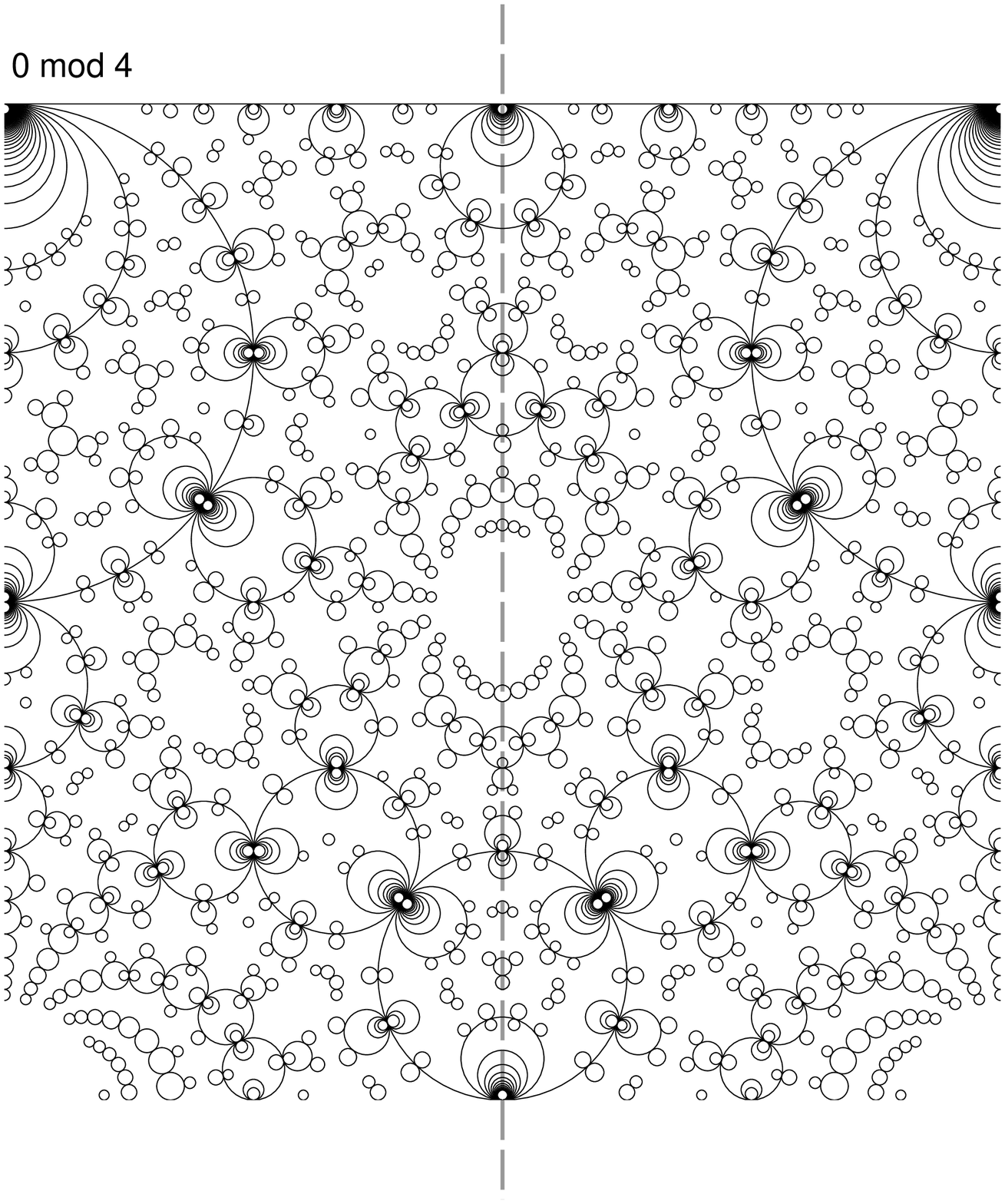}}
\caption{Circles of curvature $0~(\bmod~4)$}~\label{fig11}
\end{figure}

The figures empirically
indicate that the following extra reflection symmetries occur;
the line of symmetry is indicated on the figure  with
a dotted line.

(a) $1~(\bmod~2)$ circles are symmetrical under reflection
in the line $x= 1-y$.

(b) $2~(\bmod~4)$ circles are symmetrical under reflection
in the horizontal line $y= \frac{1}{2}$.

(c) $0~(\bmod~4)$ circles are symmetrical under reflection
in the vertical line $x = \frac{1}{2}$.

Based on this experimental evidence,
one  of the authors (CLM) conjectured that
these symmetries hold.
They  were subsequently proved by
S. Northshield ~\cite{No03}.

We may also illustrate these symmetry properties at the
level of circles of a fixed curvature. Recall from \S3
that each such circle contains a unique Apollonian packing
having it as outer circle. Such a circle
can then be  labelled by  root quadruple in the
sense of \cite{GLMWY2} of this integral Apollonian packing.  
The cases of curvature $-6, -8,$ and $-9$ corresponding to
the three cases above were pictured in \S2.
In our paper considering number-theoretic
properties of integral Apollonian packings, it was shown
in  \cite[Theorem 4.2]{GLMWY2} that the number of distinct primitive
integral Apollonian packings with a given curvature $-n$ 
of the outer circle had an interpretation as a class number
$h^{\pm}(-4n^2)$ of positive definite binary quadratic forms of discriminant
$\Delta= -4n^2$, under $GL(2, \zz)$-equivalence. One can raise 
the  question whether there is some interpretation
of the extra symmetries (a)-(c) above in terms of the 
associated class group structure under $SL(2, \zz)$ equivalence.

%
%
%
\section{The Super-Apollonian Group has Finite Covolume}
\setcounter{equation}{0}
In this section we show that the super-Apollonian group
is of finite volume as a discrete subgroup of
the real Lie group 
$Aut(Q_{\sD}) \simeq O(3, 1)$. This 
follows from  the fact that
the integer Lorentz group $O(3, 1; \ZZ)$ has
finite covolume in $O(3, 1)$, and the following result.

%
%
%
%

\begin{theorem}~\label{th51} 
(1) The super-Apollonian group $\sA^{S}$ is a normal
subgroup of index 48 in the group $Aut(Q_{D}, \ZZ)$.
The group $Aut(Q_{D}, \ZZ)$ is generated by
the super-Apollonian group and the finite
group of order 48 generated by the $4 \times 4$
permutation matrices and $\pm \bI$.

(2) The super-Apollonian group $\sA^{S}$is a normal subgroup
of index 96 in the group
$G= \bJ_0~ O(3, 1; \ZZ) \bJ_0^{-1}$, where
$$
\bJ_0=\frac{1}{2} \left[ \begin{array}{rrrr}
1 & 1 & 1 & 1 \\
1 & 1 & -1 & -1 \\
1 & -1 & 1 & -1 \\
1 & -1 & -1 & 1 \end{array} \right].
$$
The group $G$ is generated by $Aut(Q_{D}, \ZZ)$
and the duality matrix $\bD$, and consists
of matrices with integer and half-integer entries.
\end{theorem}

The duality matrix $\bD$ is given by (\ref{matrixD}) in \S5.
This result  allows us to associate to the
super-Apollonian group the normal
subgroup 
$\tilde{\sA}^{S} := \bJ_0^{-1} \sA^{S} \bJ_0$
of the integer Lorentz group $O(3, 1; \ZZ)$, of index $96$.

The proof of Theorem \ref{th51} will be
derived in a series of four lemmas.  
Let $\Gamma$ be the group generated by adjoining to $\sA$
the elements of the finite group of order $48$ given by
$\text{Perm}_4$ and $\pm I$, and let
$\tilde \Gamma=< \Gamma, \bD>$. 
The lemmas will  prove that $\Gamma= Aut(Q_D, \zz)$
and $\tilde \Gamma=\bJ_0~ O(3, 1; \ZZ) \bJ_0^{-1}$.
Then analysis of cosets of $\sA$ in these groups 
permits determining the index of $\sA$ in these groups
and showing normality. 

To determine $\Gamma$ and $\tilde{\Gamma}$ we will show
$$
\Gamma \leq Aut(Q_D, \zz) \leq \tilde{\Gamma} = G.
$$
It is easy to show that  
$\Gamma$ is a subgroup of $\tilde \Gamma$ of index 2,
and $\bD \in \tilde \Gamma$ but $\bD \notin Aut(Q_D, \zz)$,
so we get  
$Aut(Q_D, \zz) =\Gamma =<\sA^S, \text{Perm}_4, \pm \bI>$. 
It is easy to check the
inclusion  $\Gamma \leq Aut(Q_D, \zz)$ since the generators 
$\bS_i, \bS^T_i$, $\bP_\sigma$ and $\pm \bI$ of $\Gamma$ are all in 
$Aut(Q_D, \zz)$. The inclusion $Aut(Q_D, \zz) \leq G$ is proved in
Proposition \ref{lemma52}.  The equation $G=\tilde \Gamma$ is proved 
in the next three lemmas.
In Lemma \ref{lemma53} we prove that the integer Lorentz group is exactly 
the group $Aut(L_\zz)$ of 
invertible linear transformations that leave the integral Lorentz cone 
$L_\zz$ invariant.  Lemma \ref{lemma54} states that in order to show 
that a group $\sG$ of invertible 
linear transformations of $L_\zz$ equals $Aut(L_\zz)$, 
it is enough to show that (1) the action of $G$ on $L_\zz$ is transitive, 
and (2) there exists a point $v \in L_\zz$ such that the stabilizer 
$\sS_v =\{ \bU \in Aut(L_\zz) \ | \ \bU v =v \}$ is a subset of $\sG$. 
Using Lemma \ref{lemma54}, we check that (1) the action of 
$\bJ_0 \tilde \Gamma \bJ$ on $L_\zz$ is transitive, and (2), 
$\bJ_0 \tilde \Gamma \bJ$ contains the stabilizer $\sS_v$ of the point
$v =(1,1,0,0) \in L_\zz$. This proves $G= \tilde \Gamma$. 

%
%
%

\begin{lemma} \label{lemma52} It is true that
$$Aut(Q_D, \zz) \leq G= \bJ_0 O(3,1;\zz) \bJ_0^{-1}.$$
\end{lemma}
\pf
Let $\bU \in Aut(Q_D, \zz)$. We need to show that the entries of
$\bJ_0^{-1} \bU \bJ_0 =\bJ_0 \bU \bJ_0$ are all integers. Since
$\bJ_0= \frac{1}{2}\bo \bo^T -\bT$, where
\[
\bT=\left[ \begin{array}{cccc}
0 & 0 & 0 & 0 \\
0 & 0 & 1 & 1 \\
0 & 1 & 0 & 1 \\
0 & 1 & 1 & 0 \end{array} \right],
\]
we have
\beq
\bJ_0 \bU \bJ_0 = \frac{1}{4}(\sum_{i,j} U_{ij})\bo \bo^T -\frac{1}{2}\bo \bo^T
\bU \bT -\frac{1}{2} \bT \bU \bo \bo^T+\bT \bU \bT,
\eeq
where $\bT \bU \bT$ is an integer matrix.
                                                                                
From $\bU^T \bQ_D \bU=\bQ_D$ and $ \bQ_D=\frac{1}{2}(2\bI-\bo \bo^T)$, we have
\beq ~\label{6.2-for1}
\bU ^T (2\bI-\bo \bo^T) \bU=2\bI -\bo \bo^T.
\eeq
Denote by $\bv_i$ the $i$-th column of $\bU$,
and by $size(\bv):=\bo^T \bv$ the size of a vector $\bv$.  Equating the entries
of \eqref{6.2-for1} we get
$2 \bv_i \cdot \bv_j -size(\bv_i)size(\bv_j) =2\delta_{ij}-1$.
In particular, $size(\bv_i)$, $i=1,2,3, 4$  are  odd integers. It follows that
the matrix $\frac{1}{2}\bo \bo^T \bU \bT$ is integral.

Note that $Aut(Q_D, \zz)$ is closed under transposition. This is because
\beq \label{6.2for-2}
& & \bU \in Aut(Q_D, \zz) \Longrightarrow \bU^T \bQ_D \bU=\bQ_D
 \Longrightarrow \bU^T \bQ_D \bU \bQ_D \bU^T=\bU^T  \nonumber \\
& \Longrightarrow & \bU \bQ_D \bU^T =( \bU^T \bQ_D)^{-1} \bU^T=(\bQ_D)^{-1}=\bQ_D
\Longrightarrow \bU^T \in Aut(Q_D, \zz).
\eeq
Applying the same argument of the preceding paragraph to $\bU^T$, we
then prove that the matrix $\frac{1}{2} \bT \bU \bo \bo^T$ is integral.

Again using $2 \bv_i \cdot \bv_j -size(\bv_i)size(\bv_j) =2\delta_{ij}-1$, 
summing over
$i, j =1, \dots, 4$, we get
\[
(\sum_{i,j} U_{i,j})^2=2(\bv_1+\bv_2+\bv_3+\bv_4)\cdot 
(\bv_1+\bv_2+\bv_3+\bv_4) +8.
\]
Note that $\bv_1+\bv_2+\bv_3+\bv_4=(size(\br_1), size(\br_2), size(\br_3), 
size(\br_4))^T$
where $\br_i$ is the $i$-th row of the matrix $\bU$. The sum
$\sum_i (size(\br_i))^2$
is a multiple of 4 since each $ size(\br_i)$ is odd. It follows that
$(\sum_{i,j} U_{i,j})^2$ is a multiple of 8. Since $\sum_{i,j} U_{i,j}$
is an integer, we conclude that $\sum_{i,j} U_{i,j}$ is a multiple of 4
and the matrix  $\frac{1}{4}(\sum_{i,j} U_{ij})\bo \bo^T $ is integral.
This proves Lemma  \ref{lemma52}.  ~~\bsq

Clearly $\sA^S$, $\text{Perm}_4$ and $\pm \bI$ are subgroups
of $Aut(Q_D, \zz)$.  Let $\Gamma$ be the group generated by $\sA^S$,
$\text{Perm}_4$ and $\pm \bI$, and let 
$\tilde \Gamma := \langle\Gamma, \bD \rangle$.
Then $\Gamma$ is a subgroup of $\tilde \Gamma$ of index $2.$
                                                                                
The {\em Lorentz light cone} is the set of points
$\{ (y_0, y_1, y_2, y_3)^T \in \rr^4 ~: \ -y_0^2+y_1^2+y_2^2+y_3^2 = 0\}$.
 Let $L_\zz$ be the set of integer points in the
Lorentz light cone, i,e,
$$L_\zz :=\{ (y_0, y_1, y_2, y_3)^T \in \zz^4 ~:~ \ -y_0^2+y_1^2+y_2^2+y_3^2 = 0\},$$
and let $Aut(L_{\zz})$ be the set of linear transformations that
leave $L_\zz$ invariant.
%
%
%

\begin{lemma}~\label{lemma53}
 $Aut(L_\zz) =O(3, 1, \zz)$.
\end{lemma}
\pf It is clear that $O(3,1;\zz) \subseteq  Aut(L_\zz)$. 
To show the other direction,
let $\bU \in Aut(L_\zz)$. For any integer point $\bv \in L_\zz$,
$U\bv \in L_\zz$. Therefore $(\bU \bv)^T Q_\sL (\bU \bv)=0$, i.e.,
$\bv^T (\bU^T Q_{\sL} \bU) \bv=0$.
It is easy to check that the only symmetric matrices $\bQ$ satisfying
$\bv^T \bQ \bv=0$ for all $\bv \in L_\zz$ are of the form 
$\bQ=diag[-a, a, a, a]$.
Hence $\bU^T \bQ_\sL \bU=c^2 \bQ_\sL$, where $c=det(\bU)$.

Let
\[
\bX=\left[\begin{array}{cccc}
1 & 1 & 1 & 1 \\
1 & -1 & 0 & 0 \\
0 & 0 & 1 & 0 \\
0  & 0 & 0 & 1 \end{array} \right].
\]
Note that every column vector of $\bX$ is an integer point in $L_\zz$. Let
$\bY=\bU \bX$, and $\bZ=\bU^{-1} \bX$. Since $\bU \in Aut(L_\zz)$,  
every column of $\bY$
and $\bZ$ is also an integer point in $L_\zz$.
Therefore $\det(\bY)=c \det(\bX)=-2c$,
and $c \det(\bZ)=\det(\bX)=-2$. For each $(y_0, y_1, y_2, y_3)^T \in L_\zz$,
we have $y_0^2=y_1^2+y_2^2+y_3^2$. It follows that
 either $y_0$, $y_1$, $y_2$, $y_3$ are
all even, or $y_0$ and exactly one of $y_1$, $y_2$, $y_3$ are even. 
In both cases,
$\det(\bY)$ and $\det(\bZ)$ are even. This forces $c=\pm 1$. Hence
$\bU^T \bQ_\sL \bU =\bQ_{\sL}$, i.e., $\bU \in O(3,1)$.

Since $\bU$ maps points  $(1, \pm1, 0, 0)^T$, $(1, 0, \pm1,
0)^T$, $(1, 0, 0, \pm1)^T$ to integer points, if  $(a, b, c, d)$ is a row
of $\bU$, then $a\pm b$, $a\pm c$, $a\pm d \in \zz$.
Thus there exist integers $a'$, $b'$, $c'$, $d'$ of the same parity
such that $a=\frac{a'}{2}$, $b=\frac{b'}{2}$, $c=\frac{c'}{2}$, 
$d=\frac{d'}{2}$.
However, by Equation \eqref{6.2for-2}  $\bU \bQ_\sL \bU^T =\bQ_\sL$. 
This implies
 $-a^2+b^2+c^2+d^2=\pm1$. Therefore ${b'}^2+{c'}^2+{d'}^2\equiv {a'}^2$
(mod 4). Hence $a'$, $b'$, $c'$, $d'$ must all be even, 
which means  $a$, $b$, $c$, $d$
are integers, and $\bU \in O(3,1; \zz)$. ~~\bsq
%
%
%

\begin{lemma} \label{lemma54}
 Let  $\cal G$ be a group of linear transformations that preserves
$L_\zz$. If the action is  transitive, and there exists
a point $v \in L_\zz$ such that the stabilizer
$\sS_v=\{ \bU \in Aut(L_\zz) | \bU v=v\} \subseteq {\cal G}$,
then ${\cal G}=Aut(L_\zz)$.
\end{lemma}
\pf  Clearly $\sG \subseteq Aut(L_\zz)$.
For any $\bP \in Aut(L_\zz)$, assume $\bP(v)=v'$. Since
$\sG$ acts transitively on $L_\zz$, there exists $\bG \in \sG$ such
that $\bG(v')=v$. That is,  $\bG \bP(v)=v$. So $\bG \bP 
\in \sS_v \subseteq \sG$, and
then $\bP \in \bG^{-1} \sG=\sG$. This proves $Aut(L_\zz) \subseteq
\sG$.~~~\bsq

%
%
%

\begin{lemma} \label{lemma55}
  $\tilde \Gamma = G= \bJ_0 O(3,1;\zz)\bJ_0^{-1}$.
\end{lemma} 
                
\pf
    By Lemma \ref{lemma53}, it is  sufficient to prove that
\beq \label{6.2for-3}
\bJ_0^{-1} \tilde \Gamma \bJ_0=\bJ_0 \tilde \Gamma \bJ_0 =Aut(L_\zz).
\eeq
It is straightforward to check that $\bJ_0 \bS_i \bJ_0$, $\bJ_0 \bS^T_i \bJ_0$,
$\bJ_0 \bP_{\sigma} \bJ_0$ and $\bJ_0 \bD \bJ_0$ are integer matrices. 
In particular,
\beqs
\bJ_0 \bS_1 \bJ_0=\left[ \begin{array}{rrrr}
          2 & -1 & -1 & -1 \\
          1 &  0 & -1 & -1 \\
          1 & -1 &  0 & -1 \\
          1 & -1 & -1 &  0 \end{array} \right],
\eeqs
$\bJ_0 \bP_{34} \bJ_0 =\bP_{34}$ and $\bJ_0 \bD \bJ_0 =diag[1,-1,-1,-1]$. 
Therefore
$\bJ_0 \tilde \Gamma \bJ_0 \subseteq O(3,1;\ZZ) = Aut(L_\zz)$.
 
   For any integer point $(y_0, y_1, y_2, y_3) \in L_\zz$,
$-y_0^2+y_1^2+y_2^2+y_3^2=0$ implies \\
$y_0+y_1+y_2+y_3\equiv 0~ (\bmod~ 2)$.
Hence $\bJ_0 (y_0, y_1, y_2, y_3)^T$ is integral. It follows that
\[
\bJ_0(L_\zz)=
\{(a_1, a_2, a_3, a_4) \in \zz^4 ~:~ (a_1, a_2, a_3, a_4) 
\text{ are curvatures of a Descartes configuration }\}.
\]
Then Theorem \ref{reduction}  implies that 
$\bJ_0 \tilde \Gamma \bJ_0$ acts transitively
on the integral Lorentz light  cone $L_\zz$, where $-\bI$ exchanges the total
orientation of a point in $L_\zz$. 
                                                                              
We use Lemma \ref{lemma54} to prove equation \eqref{6.2for-3}
with $\sG=J_0 \tilde \Gamma J_0$.
Let $\bv :=(1, 1, 0,0) \in L_\zz$ and consider its stabilizer
$$\sS_{\bv} :=\{ \bU \in O(3, 1;\zz) ~:~ \bU \bv=\bv, \ \bU^T \bQ_\sL \bU 
= \bQ_\sL\}.$$
Assume $\bU=(u_{i,j})_{i,j=1}^4$.
Solving the equations
      \beq
\bU(1,1,0,0)^T=(1,1,0,0), \qquad
\bU^T \bQ_\sL \bU= \bQ_\sL,
\eeq
we obtain the following linear and quadratic relations
between the entries of $\bU$:
\beqs
u_{12}=1-u_{11},~~~~~ u_{13}&= & u_{23},~~~~~ u_{14}=u_{24}, \\
u_{22}=1-u_{21},~~~~ u_{32} &= & -u_{31},~~~~ u_{42}=-u_{41},
\eeqs
and
\beqs
u_{33}^2+u_{43}^2=u_{34}^2+u_{44}^2=1,& &\ u_{33}u_{34}+u_{43}u_{44}=0, \\
u_{13}=u_{31}u_{33}+u_{41}u_{43},  & &\ u_{14}=u_{31}u_{34}-u_{41}u_{44}, \\
u_{21}=\frac{1}{2}(u_{31}^2+u_{41}^2).~~~
\eeqs
  It follows that the matrix $\bU$ can be expressed as
\beq\label{form}
\bU=\left( \begin{array}{cccc}
1+t & -t & gm+hn & km+ln \\
t & 1-t & gm+hn & km+ln \\
m & -m & g & k \\
n & -n & h & l \end{array}\right),
\eeq
where $t=(m^2+n^2)/2$, $g^2+h^2=k^2+l^2=1$, and $gk+hl=0$.
Since $g, h, k, l \in \zz$, we must have
 $(g, h), (k, l) \in \{(\pm 1, 0), (0, \pm 1)
\}$.
                                   
We can classify the matrices of the form \eqref{form} into four types,
up to a possible  multiplication by $\bP_{34}$)
, as follows.
\beqs
\text{Type I}: \ \left( \begin{array}{cccc}
 1+t & -t & m & n\\
 t & 1-t & m & n\\
 m & -m & 1 & 0\\
 n & -n & 0 & 1 \end{array}\right),
\qquad
\text{Type I\@I}: \ \left( \begin{array}{cccc}
 1+t & -t & m & -n\\
 t & 1-t & m & -n\\
 m & -m & 1 & 0\\
 n & -n & 0 & -1 \end{array}\right),
\eeqs
\beqs
\text{Type  I\@I\@I}: \ \left( \begin{array}{cccc}
 1+t & -t & -m & n\\
 t & 1-t & -m & n\\
 m & -m & -1 & 0\\
 n & -n & 0 & 1 \end{array}\right),
\qquad
\text{Type I\@V}: \ \left( \begin{array}{cccc}
 1+t & -t & -m & -n\\
 t & 1-t & -m & -n\\
 m & -m & -1 & 0\\
 n & -n & 0 & -1 \end{array}\right),
\eeqs
where $t=(m^2+n^2)/2$ and $m,n, t \in \zz$.
                                                                                
We denote  a matrix of type $X$ with parameters $m,n$ by $\bU(m,n; {\rm X})$.
The following equations can be easily checked.
\beqs
 \bU (m,n; {\rm I}) \bU (k,l; {\rm I})&=& \bU (m+k, n+l; {\rm I}),\\
 \bU (m,n;\text{I\@I})  \bU (k,l; {\rm I}) & = &
         \bU (m+k, n-l; \text{I\@I}),\\
 \bU (m,n; \text{I\@I\@I})  \bU (k,l; \text{I}) & = &
        \bU (m-k, n+l; \text{I\@I\@I}),\\
 \bU (m,n; \text{I\@V})  \bU (k,l; \text{I}) & = &
         \bU (m-k, n-l; \text{I\@V}).
\eeqs
Also we have
\[
 \bU (1,1; {\rm I})^{-1} =\bU (-1,-1;{\rm I}), \qquad
 \bU (1, -1; {\rm I})^{-1}= \bU (-1, 1; {\rm I}).
\]
Therefore the stabilizer $\sS_v$ is generated by $\bP_{34}$,
$\bU (0,0; {\rm I})$, $\bU (0,0; \text{I\@I})$
$\bU (0,0; \text{I\@I\@I})$, $\bU (0,0; \text{I\@V})$
together with $\bA= \bU (1,1;{\rm I}), \bB=\bU (1, -1;{\rm  I})$.
                                                                                
Note that $\bU (0,0; \text{I})=\bU (0,0; \text{I\@I})^2$,
$ \bU (0,0; \text{I\@I\@I})=\bP_{34} \bU (0,0; \text{I\@I}) \bP_{34}$,
and $\bU (0,0; \text{I\@V})=\bU (0,0; \text{I\@I})
\bU (0,0; \text{I\@I\@I})$.
Thus $\sS_v$ is generated by the matrices $\bP_{34}$, $\bA$, $\bB$, and
$\bC=\bU (0,0; \text{I\@I})$.

Now  one can check that
\beqs
\bP_{34} & = & \bJ_0 \bP_{34} \bJ_0 ~\in \bJ_0 \tilde \Gamma \bJ_0, \\
\bC & = & {\rm diag}[1,1,1,-1] =\bJ_0(-\bI \bP_{23} \bP_{14} \bD)\bJ_0 
~\in \bJ_0 \tilde \Gamma \bJ_0,  \\
\bA & = & (\bJ_0 \bS_1 \bJ_0 ){\rm diag}[1,1,-1, -1]= 
\bJ_0 (\bS_1 \bP_{12} \bP_{34}) \bJ_0
~\in \bJ_0 \tilde \Gamma \bJ_0,\\
\bB & = & \big(\bJ_0 \bP_{34}  \bJ_0\big)  \bC \big(\bJ_0 \bS_2 \bJ_0  \big) \bC ~\in
 \bJ_0 \tilde \Gamma \bJ_0.
\eeqs
Hence $\sS_v \subseteq  \bJ_0 \tilde \Gamma \bJ_0.$ By Lemma \ref{lemma54}
$ \bJ_0 \tilde \Gamma \bJ_0=Aut(L_\zz)=O(3,1;\zz)$, or equivalently,
$\tilde \Gamma =\bJ_0 O(3,1;\zz) \bJ_0^{-1}$. This finishes the proof.
~~~\bsq \\
%
%
%

\vspace{.3cm}
\noindent {\bf Proof of Theorem~\ref{th51}.}
Lemma \ref{lemma52} and Lemma \ref{lemma55} show that
$\Gamma \leq Aut(Q_D, \zz) \leq G$. Hence $\Gamma =Aut(Q_D, \zz)$ since
$\Gamma$ is a subgroup of $G$ of index 2 and
the duality matrix $\bD$ is not in $Aut(Q_D, \zz)$. In other words,
$Aut(Q_D, \zz)$ is generated by the super-Apollonian group $\sA^S$ and the
finite group of order 48 generated by  the $4 \times 4$ permutation matrices
and $\pm I$.
The super-Apollonian group $\sA^S$ is a normal subgroup
of $Aut(Q_D, \zz)$,
since $\bP_\sigma \bS_i \bP_{\sigma^{-1}} =\bS_{\sigma(i)}$, and
$\bP_\sigma \bS_i^T \bP_{\sigma^{-1}} =\bS_{\sigma(i)}^T$. The index is 48
since $\sA^S \cap (\text{Perm}_4 \times \{ \pm \bI\})=\bI$, (c.f. \S5
of Part I \cite{GLMWY11}).

By Lemma \ref{lemma55}, the group $G$ is generated by $Aut(Q_D, \zz)$
and $\bD$. Note that $\bD^2=I$ and $\bD \bS_i \bD =\bS^T_i$. 
It follows that the super-Apollonian group is a normal  subgroup 
of $G$ with index 96. ~~\bsq \\

The second part of Theorem~\ref{th51} can be rephrased
as asserting that the conjugate group
$\tilde{\sA}^{S} = \bJ_0 \sA \bJ_0$ is a normal subgroup of
index $96$ in $O(3,1)$. Its generators 
$\tilde{\bS}_j = \bJ_0 \bS_j \bJ_0^{-1}$
and $\tilde{\bS}_j^{\perp} = \bJ_0 \bS_j^{\perp} \bJ_0^{-1}$
are given by
$$
\tilde{\bS}_1 
\left[ \begin{array}{rrrr}
2 & -1 & -1 &-1 \\
1 & 0 & -1 & -1 \\
1 & -1 & 0 & -1 \\
1 & -1 & -1 & 0
\end{array} 
\right] ~~~\mbox{and}~~
\tilde{\bS}_2  = \left[ \begin{array}{rrrr}
2 & -1 & 1 &1 \\
1 & 0 & 1 & 1 \\
-1 & 1 & 0 & -1 \\
-1 & 1 & -1 & 0
\end{array}
\right] 
$$
and
$$
\tilde{\bS}_3 = \left[ \begin{array}{rrrr}
2 & 1 & -1 &1 \\
-1 & 0 & 1 & -1 \\
1 & 1 & 0 & 1 \\
-1 & -1 & 1 & 0
\end{array} 
\right] ~~~\mbox{and}~~
\tilde{\bS}_4  = \left[ \begin{array}{rrrr}
2 & 1 & 1 &-1 \\
-1 & 0 & -1 & 1 \\
-1 & -1 & 0 & 1 \\
1 & 1 & 1 & 0
\end{array}
\right] 
$$
The generators $\tilde{\bS}_j^{\perp} = \bJ_0 \bS_j^{\perp} \bJ_0^{-1}$
are given by 
$$
\tilde{\bS}_j^{\perp} = (\tilde{\bS}_j)^{T},
$$
which follows using $\bS_j^{\perp} = \bS_j^T$ 
and $\bJ_0 = \bJ_0^T = \bJ_0^{-1}.$

%
%
%
\section{Super-Integral Super-Packings}
\setcounter{equation}{0}

This section  treats the strongest form of integrality for
super-packings, which is that where one
(and hence all) Descartes configurations in the
super-packing have 
an integral augmented curvature-center coordinate
matrix $\bW_{\sD}$.
We say that a Descartes configuration with this
property is {\em super-integral}, and the
same for the induced super-Apollonian packing.
The following result classifies such packings.
%
%
%
%

\begin{theorem}~\label{th81n}
(1) These are  exactly $14$ different geometric
super-packings that are super-integral.

(2) The set of ordered, oriented Descartes
configurations that are super-integral comprise
$672$ orbits of the super-Apollonian
group.
\end{theorem}

These packings are classified  here as rigid objects, not movable
by Euclidean motions. 
To prove this result, it suffices to determine which strongly-integral
configurations are super-integral.
The next result classifies 
the possible types of super-integral Descartes configurations,
according to the allowed value of their divisors.

%
%
%
%

\begin{theorem}~\label{super-int}
Suppose that
an ordered, oriented  Descartes configuration $\sD$ in $\RR^2$ has integral 
curvature-center
coordinates $\bM=\bM_{\sD}$, and let $g= \gcd(a_1, a_2, a_3, a_4)$,
where $(a_1, a_2, a_3, a_4)$ is its first column of signed
curvatures.
Then $\sD$ has integral augmented curvature-center coordinates
$\bW_{\sD}$ if and only if one of the following conditions hold. 

(i) $g=1,$ or

(ii) $g= 2$, and each row of $\bM$ has 
the sum of its last two entries being $1~(\bmod~2)$. 

(iii) $g=4$, and the last columns rows of $\bM$ are congruent to
$$
\left[ \begin{array}{cc}
1 & 0 \\
1 & 0 \\
1 & 0 \\
1 & 0 \\ \end{array} \right]  ~~~\mbox{or}~~~
\left[ \begin{array}{cc}
0 & 1 \\
0 & 1 \\
0 & 1 \\
0 & 1 \\ \end{array} \right] ~~(\bmod~2).
$$
\end{theorem}
\pf
By Theorem \ref{normal_form}, there exists a matrix $\bU \in \sA^S$ and 
a permutation matrix $\bP$ such that 
\[
\bP \bU \bM_\sD =A_{m,n}[g] \text{ or } B_{m,n}[g],
\]
for some $m, n \in \{0,1\}$, 
where $A_{m,n}[g]$ and $B_{m,n}[g]$ are given in Formula \eqref{ABmn}, 
and their
corresponding augmented curvature-center coordinate matrices are
$\tilde A_{m,n}[g]$, $\tilde B_{m,n}[g]$, given in Formula \eqref{aug-ABmn}. 
Since each generator of $\sA^S$ preserves the super-integrality, as
well as the parity of every element of $\bW_\sD$,  it follows that 
$\bW_\sD$ is integral if and only if one of the following conditions
holds:
\begin{enumerate}
\item $g=1$, or
\item $g=2$, and $m^2+n^2-1 \equiv 0~(\bmod~2)$, or
\item $g=4$, and the reduced form is $A_{0,1}[4]$ or $B_{1,0}[4]$, up 
to a permutation of rows.  
\end{enumerate}
In Case 2 we need $m+n \equiv 1 (\bmod 2)$; in Case 3 the condition is
equivalent to the one stated in the theorem. 
~~~$\bsq$

\paragraph{Proof of Theorem~\ref{th81n}.} 
(1) Since each geometric super-packing corresponds to 48 distinct orbits of
the super-Apollonian group on ordered, oriented Descartes configurations, 
to show  there are exactly 14 geometric super-packings,
it suffices to show there are exactly 672 orbits of the super-Apollonian
group that 
are super-integral, which is (2).

(2)Theorems \ref{normal_form} and \ref{super-int} allow us to classify the 
set of ordered, oriented Descartes configurations that are super-integral 
by the action of $\sA^S$.  
From the criterion of Theorem \ref{super-int}, we have:
\begin{enumerate}
\item For $g=1$, any strongly integral Descartes configuration is 
super-integral. 
\item For $g=2$, half of those orbits are super-integral, namely, those
 whose reduced forms are $A_{0,1}[2], A_{1,0}[2], B_{0,1}[2]$ or $B_{1,0}[2]$, 
up to 
 a permutation of rows.
\item For $g=4$, one fourth of those orbits are super-integral, namely, 
those whose reduced form are $A_{0,1}[4]$ or $B_{1,0}[4]$, up to a permutation 
of rows. 
\item For $g \neq 1, 2, 4$, there are no super-integral Descartes 
configurations.   
\end{enumerate}

\begin{table}[htpb] 
\begin{center}
\begin{tabular}{|c|c|c|} 
\hline 
$\bg$  & \# of orbits of $\sA^S$  & Representative \\
\hline
$(1,1,1,1)$   & 96   & $A_{1,1}[1],~~~ B_{1,1}[1]$ \\
$(2,1,1,1)$   & 96   & $A_{1,0}[1], ~~~B_{0,1}[1]$ \\
$(1,1,2,1)$   & 48   & $A_{0,0}[1]$              \\ 
$(1,1,1,2)$   & 48   & $B_{0,0}[1]$  \\
$(4,1,2,1)$    & 48   & $A_{0,1}[1]$  \\
$(4,1,1,2)$    & 48   & $B_{1,0}[1]$  \\ \hline 
$(1,2,1,1)$    & 96   & $A_{1,0}[2], ~~~B_{0,1}[2]$ \\
$(2,2,2,1)$    & 48    & $A_{0,1}[2]$  \\
$(2,2,1,2)$    & 48     & $B_{1,0}[2]$  \\ \hline
$(1,4,2,1)$    & 48     & $A_{0,1}[4]$  \\
$(1,4,1,2)$    & 48     &  $B_{1,0}[4]$ \\
\hline
\end{tabular}
\caption{Orbits of super-integral Descartes configurations classified by $\bg$.}
~\label{ta81}
\end{center}
\end{table} 

More details of this calculation are given in Table~\ref{ta81}.
To explain the notation in Table~\ref{ta81}, 
for any $4\times 4$ integral matrix $\bW$, let $g_i$ be the greatest common 
divisor of entries $w_{1,i}$, $w_{2,i}$, $w_{3,i}$, $w_{4,i}$ 
in the $i$-th column. 
Then the action of $\sA^S$ preserves the vector $\bg=(g_1, g_2, g_3, g_4)$. 
(For $\bW=\bW_\sD$, $g_2$ is the greatest common divisor of the curvatures.)  
In each  row of the table  we give the number of orbits of $\sA^S$ formed by 
the set of ordered, oriented  Descartes configurations that are super-integral
with the given $\bg$. We also list the representatives of those orbits. 
Note that each entry in the column labelled ``Representative'' 
stands for 48 orbits, obtained 
by taking two choices of (total) orientation, and 24 choices of permutation
of rows. Also note the symmetry that the 
configurations with $\bg$  are inverse of the ones with
$\bg'=\bP_{12} \bg$.
We  conclude from the table that the set of  ordered, oriented
Descartes configurations that are super-integral comprise  
$384(1+\frac{1}{2} +\frac{1}{4})=672$ orbits of the super-Apollonian group.
~~~$\bsq$

%
%

\section{Concluding remarks}
\setcounter{equation}{0}

This paper showed that the ensemble
of all primitive,
strongly integral Apollonian circle packings can  be
simultaneously described 
in terms of  an orbit of a larger discrete group,
the super-Apollonian group, acting on the
standard strongly-integral super-packing.  
Study of the  locations of
the individual integer packings inside 
the standard  super-packing, presented in \S6,
 leads  to interesting
questions, not all of which are resolved.
The standard super-packing also played a
role in analyzing the structure of the
super-Apollonian group as a discrete subgroup,
carried out in  \S7.

The various illustrations show 
 the usefulness of graphical representations, 
as a guide to both finding and illustrating
results. This contribution  is due in large part
to the statistician co-authors (CLM and AW).
Graphics were  particularly useful in finding extra
symmetries of these objects, such as those illustrated
in Figures 9--11 and subsequently proved by
S. Northshield \cite{No03}.
However one must not forget the
adage of H. M. Stark \cite[p. 225]{St78}:
``A picture is worth a thousand words, provided 
one uses another thousand words to justify
the picture.''  Section 3 of this paper 
provides such a justification for certain
features of Figure 4.

There remain some open questions, particularly
concerning the classification of all 
integer root quadruples 
 classified by fixed curvature $-N$ of the
bounding circle. This quantity is known
to be interpretable as a class number,
as described in \cite[Theorem 4.2]{GLMWY2}.
In \S6 we observed some symmetries of these
root quadruples inside the standard super-packing,
see Figures ~\ref{fig4}--\ref{fig6}.  
There is a new invariant
that can be associated to such quadruples,
which is their nesting depth as defined in \S4
with respect to the generating quadruple
$\sD_0$ of the standard super-packing. It would
be interesting to see whether this invariant
might give some  further insight into class numbers.

%
%
\newpage
\section{Appendix.  Strong Integrality Criterion}
\setcounter{equation}{0}

This appendix  establishes that to
show a Descartes configuration is strongly integral
it suffices to show that
 three of its four circles are strongly integral.
This affirmatively answers a question posed to us by K. Stephenson.

%
%
%
%

\begin{theorem}~\label{th101}
A Descartes configuration $\sD$ has integral curvature-center
coordinates $\bM_{\sD}$ if and only if it contains three
circles having integer curvatures and whose
curvature$\times$centers, viewed as complex numbers,  lie in $\zz[i]$.
\end{theorem}
                                                                                
\pf
The condition is clearly necessary, and the problem is
to show it is sufficient. We write the circle
centers as complex numbers $\bz_j = x_j + i y_j$. So suppose $\sD$
contains three circles with curvatures
$b_1,~b_2, ~b_3 \in \zz$ and with
curvature$\times$centers $b_1\bz_1,~b_2\bz_2,~b_3\bz_3 \in \zz[i]$.
We must show that
the fourth circle in the configuration has
$b_4 \in \zz$ and $b_4\bz_4 \in \zz[i].$

For later use, we note that Theorem 3.1 of Part I \cite{GLMWY11} has a nice
interpretation using complex numbers $\bz$
to represent circle centers. This was formulated  in \cite{LMW02}
as the Complex Descartes theorem. 
It gives 
\beql{901}
b_1^2\bz_1^2 +  b_2^2\bz_2^2 + b_3^2\bz_3^2 +  b_4^2\bz_4^2 =
\frac{1}{2}( b_1\bz_1 +  b_2\bz_2 + b_3\bz_3 +  b_4\bz_4)^2.
\eeq
and
\beql{soddy}
b_1^2\bz_1 +  b_2^2\bz_2 + b_3^2\bz_3 +  b_4^2\bz_4
= \frac{1}{2} ( b_1\bz_1 +  b_2\bz_2 + b_3\bz_3 +  b_4\bz_4)
(b_1 + b_2 +b_3 +b_4).
\eeq
                                                                                
  We claim  that $b_4 \in \zz.$ This is proved in
the following two cases.

{\em Case 1. $b_1b_2b_3\neq 0$.}
                                                                                
We first suppose
that  $\bz_1=0$.
If both $x_2$ and $x_3$ are zero, then $-\frac{1}{b_1} =\frac{1}{b_2}+\frac{1}{b_3}$, which means $b_1b_2+b_2b_3+b_3b_1=0$. Hence $b_4=b_1+b_2+b_3 \in \zz$.
Otherwise by permuting $\bz_2$ and $\bz_3$ if
necessary we may assume that $x_2 \ne 0$.
Then the following equations encode the distance between 
the circle centers, since the circles touch.
\begin{eqnarray}
x_2^2+y_2^2=(\frac{1}{b_1} + \frac{1}{b_2})^2, \label{xy1_sqr} \\
x_3^2+y_3^2=(\frac{1}{b_1} + \frac{1}{b_3})^2, \label{xy_sqr} \\
(x_3-x_2)^2+(y_3 - y_2)^2
=(\frac{1}{b_2}+\frac{1}{b_3})^2. \label{xyxy1}
\end{eqnarray}
We wish to solve these equations for $y_3$ in terms of
$b_1, b_2, b_3, x_1$ and $x_2$.
To this end we subtract the
first two equations from the third and obtain
$$
2(x_2x_3 + y_2 y_3) = (\frac{1}{b_1} + \frac{1}{b_2})^2
+ (\frac{1}{b_1} + \frac{1}{b_3})^2 -(\frac{1}{b_2}+\frac{1}{b_3})^2 := R.
$$
Calling the right side of this equation $R$, we obtain
$$
x_3 = \frac{1}{2x_2}\left( R - 2y_2 y_3 \right).
$$
where the division  is allowed since $x_2 \ne 0$.
Substituting this in the second equation yields a quadratic
equation in $y_3$, with $x_3$ eliminated, namely
(after multiplying by $4x_2^2$),
$$
4(x_2^2 + y_2^2)  y_3^2 - 4R~ y_2y_3+R^2
- 4x_2^2  (\frac{1}{b_1} + \frac{1}{b_3})^2 = 0.
$$
Since $y_3$ is
rational, this  equation has rational  solutions, so the discriminant $\Delta$
must be the square of a rational number. After some
calculation one obtains
\[
\Delta=16^2 x_2^2 \left(\frac{ b_1b_2+b_1b_3+b_2b_3}{ (b_1b_2b_3)^2} \right).
\]
Since $x_2, b_1, b_2$ are
nonzero it  follows that  $b_1b_2+b_2b_3+b_1b_3$ is a perfect square.
Viewing  the Descartes equation as a quadratic equation in $b_4$
we obtain the formula
$b_4 :=b_1+b_2+b_3\pm2\sqrt{b_1b_2+b_2b_3+b_1b_3},$
which shows that both roots are integers.
These are the oriented curvatures of the
two possible choices for the fourth circle in the
Descartes configuration, so  $b_4 \in \ZZ$.
                                                                                
Now assume that   $\bz_1= x_1 + i y_1$ is arbitrary. Define
\begin{equation*}
(s_2, t_2) := (x_2 - x_1, y_2 - y_1) ~~\mbox{and}~~ (s_3, t_3) := (x_3 - x_1,
y_3-y_1).
\end{equation*}
Then $s_2, t_2, s_3, t_3$ are rational numbers, and they also satisfy
equations \eqref{xy1_sqr}, \eqref{xy_sqr}, \eqref{xyxy1}. (just replace
$x_2, y_2, x_3, y_3$ in \eqref{xy1_sqr}, \eqref{xy_sqr},
\eqref{xyxy1} by $s_2,t_2,s_3, t_3$, respectively.)
By the preceding argument, we again have
$16s_2^2 (b_1b_2+b_2b_3+b_1b_3)/(b_1b_2b_3)^2 = q^2$
for some rational number $q$.
This implies that $b_1b_2+b_2b_3+b_1b_3$ is a perfect square and
consequently that
$b_4$ is an integer.
                                                                                
{\em Case 2. $b_1b_2b_3=0$.}
This  is proved similarly,  with an
easier calculation.
                                                                                
The claim follows, so $b_4 \in \ZZ$.

We now proceed to  show that
$b_4\bz_4 \in \zz[i]$. Now  \eqn{901}  gives
\begin{equation}~\label{reduce0}
b_4\bz_4=b_1\bz_1+b_2\bz_2+
b_3\bz_3\pm2\sqrt{b_1\bz_1b_2\bz_2+b_2\bz_2b_3\bz_3+b_1\bz_1b_3\bz_3}.
\end{equation}
The equation  \eqn{soddy} gives
\begin{equation}~\label{reduce}
(b_1+b_2+b_3-b_4)b_4\bz_4=2(b_1^2\bz_1+b_2^2\bz_2+b_3^2\bz_3)-(b_1+b_2+b_3+b_4)
(b_1\bz_1+b_2\bz_2+b_3\bz_3).
\end{equation}
We treat two mutually exhaustive cases.
                                                                                
{\em Case 1.  $b_1+b_2+b_3 \neq b_4.$}
                                                                                
Then \eqn{reduce} gives
$b_4\bz_4=x_4+iy_4$ for some rational
numbers $x_4, y_4$. However
  $b_1\bz_1, b_2\bz_2, b_3\bz_3$ are integers by hypothesis,
whence \eqn{reduce0} shows that $b_4z_4$ is an algebraic integer.
Since $x_4, y_4$  are rational, we conclude that
$b_4\bz_4$ must be an integer, i.e., in $\zz[i]$.
                                                                                
{\em Case 2.  $b_1+b_2+b_3=b_4$.}
                                                                                In this case, we have
$b_1b_2+b_2b_3+b_1b_3 =0$. Now from
\eqref{reduce}, we have
\[
 b_1^2\bz_1+b_2^2\bz_2+b_3^2\bz_3=(b_1+b_2+b_3)(b_1\bz_1+b_2\bz_2+b_3\bz_3),
\]
which can be simplified to
\[
b_2b_3\bz_1+b_1b_3\bz_2+b_1b_2\bz_3=0.
\]
Thus
\begin{eqnarray*}
(b_1\bz_1)(b_2\bz_2)+(b_2\bz_2)(b_3\bz_3)+(b_3\bz_3)(b_1\bz_1) & = &
b_1(b_2\bz_2+b_3\bz_3)\frac{b_1b_3\bz_2+b_1b_2\bz_3}{-b_2b_3}
+(b_2\bz_2)(b_3\bz_3) \\
&=&-b_1^2\bz_2^2-b_1^2\bz_3^2+\bz_2\bz_3 \left(b_2b_3-
\frac{b_1^2}{b_2b_3}(b_2^2+b_3^2)\right) \\
&=& -b_1^2(\bz_2-\bz_3)^2.
\end{eqnarray*}
Now $b_1(\bz_2-\bz_3)$ is a Gaussian
rational number whose square is integral.
Hence $b_1(\bz_2-\bz_3)$ must be integral.
It follows that $b_4\bz_4=b_1\bz_1+b_2\bz_2+b_3\bz_3\pm 2b_1(\bz_2-\bz_3)i$ is
in $\zz[i]$.~~~$\bsq$ \\
                                                                                
Since the Apollonian group consists of integer matrices,
all Descartes configurations
in an  Apollonian packing generated by  a strongly integral
Descartes configuration are  strongly
integral.  This explains
the integrality properties of curvatures and
curvature $\times$ center pictured in the packing
in \S1, for example. The previous result now  gives
a weaker necessary and sufficient condition for an Apollonian
packing to be strongly integral.

%
%

\begin{theorem} \label{th32} 
An Apollonian circle packing is strongly integral if
and only if it contains three mutually tangent
circles which have integer curvature-center
coordinates.
\end{theorem}
                                                                                
\pf
Suppose we are give three mutually tangent circles
in the packing that are strongly integral.
Any set of three mutually tangent circles in
the packing is part of some Descartes configuration
in this packing. This follows from the recursive
construction of the packing, which has a finite
number of circles at each iteration. If iteration $j$
is the first iteration at which all three tangent
circles are present, at that iteration they
necessarily belong to a unique Descartes configuration.
Theorem~\ref{th101} now implies that this
Descartes configuration is strongly integral.
It then follows that the whole Apollonian
packing is strongly integral. This proves the
``if'' direction, and the ``only if'' direction
is immediate.
~~~$\bsq$

%
%
%

{\tt
\begin{tabular}{lllll}
email: & graham@ucsd.edu \\
& lagarias@umich.edu \\
  & colinm@research.avayalabs.com \\
 & allan@research.att.com \\
& cyan@math.tamu.edu
\end{tabular}
 }

\end{document}